\documentclass{article}

\usepackage{arxiv}
\usepackage[colorlinks=true, linkcolor=blue, citecolor=black, urlcolor=blue]{hyperref}
\usepackage[utf8]{inputenc} 
\usepackage[T1]{fontenc}    
\usepackage{url}            
\usepackage{booktabs}       
\usepackage{amsfonts}       
\usepackage{nicefrac}       
\usepackage{microtype}      
\usepackage{lipsum}		
\usepackage{graphicx}
\usepackage{doi}

\RequirePackage{amsthm,amsmath,amsfonts,amssymb,comment,cancel}
\RequirePackage{graphicx}
\usepackage{cleveref}
\usepackage{enumitem}

\allowdisplaybreaks

\theoremstyle{plain}

\newtheorem{theorem}{Theorem}[section]
\newtheorem{lemma}[theorem]{Lemma}
\newtheorem{proposition}[theorem]{Proposition}
\newtheorem{corollary}[theorem]{Corollary}
\theoremstyle{definition}

\theoremstyle{remark}

\providecommand{\supp}[1]{\operatorname{supp}(#1)}
\providecommand{\sw}{\operatorname{SW}}
\providecommand{\w}{\operatorname{W}}
\providecommand{\sg}{\operatorname{sgn}}

\title{An improved central limit theorem for the empirical sliced Wasserstein distance}

\author{ \href{https://orcid.org/0000-0002-8572-2541}{\includegraphics[scale=0.06]{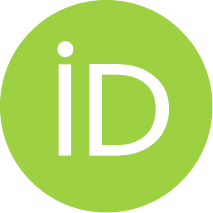}\hspace{1mm}David Rodríguez-Vítores}\thanks{Correspondence to: \texttt{david.rodriguez.vitores@uva.es}} \\
	Universidad de Valladolid and IMUVa\\
	Valladolid, Spain\\
	\And
    \href{https://orcid.org/0000-0003-3764-5411}{\includegraphics[scale=0.06]{orcid.pdf}\hspace{1mm}Eustasio del Barrio Tellado} \\
	Universidad de Valladolid and IMUVa\\
	Valladolid, Spain\\
    \And
	\href{https://orcid.org/0000-0002-1252-2960}{\includegraphics[scale=0.06]{orcid.pdf}\hspace{1mm} Jean-Michel Loubes} \\
	Institut de Mathématiques de Toulouse, INRIA\\
    Toulouse, France
}

\date{}

\renewcommand{\headeright}{}
\renewcommand{\undertitle}{}
\renewcommand{\shorttitle}{An improved CLT for the empirical sliced Wasserstein}

\begin{document}
\maketitle

\begin{abstract}
Wasserstein distances are widely used in modern data analysis but pose significant computational and statistical challenges in high dimensions. The sliced Wasserstein distance alleviates these challenges by leveraging one-dimensional projections. Building on the Efron-Stein inequality---a technique proven effective in related problems---and a non-trivial control of the optimal transport potentials across directions, we establish a central limit theorem for the $p$-sliced Wasserstein distance, for $p>1$, centered at the expected empirical cost. Unlike for the general Wasserstein distance, the centering can be replaced by the population cost, enabling valid statistical inference. This generalizes and refines existing one-dimensional results, providing the first asymptotically valid inference framework for the sliced Wasserstein distance between possibly non-compact measures. Finally, we address other practical aspects crucial for inference, including Monte Carlo approximation of the slicing integral and consistent variance estimation.
\end{abstract}


\section{Introduction}

In recent years, the rapid growth in both the quantity and variety of data has driven the development of new methods to tackle emerging challenges, highlighting the need for effective tools to analyze distributional data. Optimal transport has emerged as a powerful approach in numerous applications, including domain adaptation \cite{chang2022unified}, fair learning \cite{deLara.et.al.JMLR.2024,gordaliza2019obtaining}, generative modeling \cite{arjovsky2017wasserstein,pooladian2024pluginestimationschrodingerbridges}, representation learning \cite{tolstikhin2018wasserstein}, computer vision \cite{BonneelDigne2023}, and flow cytometry gating \cite{del2020optimalflow,freulon2023cytopt}. Formally, given $P,Q \in \mathcal P(\mathbb R^d)$, the set of probability measures on $\mathbb R^d$ with the Borel $\sigma$-field, and a cost function $c:\mathbb R^d \times \mathbb R^d \to [0,\infty)$, the optimal transport problem consists in solving
\begin{equation}\label{eq:OT}
    \min_{\pi \in \Pi(P,Q)} \int_{\mathbb R^d \times \mathbb R^d} c(x,y) d\pi(x,y) \ ,
\end{equation}
where $\Pi(P,Q)$ represents the set of probability measures with fixed marginals $P$ and $Q$. A minimizer $\pi$ in \eqref{eq:OT} is known as an \textit{optimal transport plan} or \textit{optimal transport coupling}. Intuitively, $\pi(x,y)$ represents the amount of mass transported from $x$ to $y$ in the optimal solution. If there exists a map $T:\mathbb R^d \rightarrow \mathbb R^d$ satisfying $\pi = (I_d,T)\# P$, where $I_d$ is the identity map and $\#$ denotes the push-forward operator, then $T$ is called an \textit{optimal transport map}. For the particular choice $c(x, y) = \| x - y \|^p$, with $p \geq 1$, the minimal transport cost defines the so-called \textit{$p$-Wasserstein distance} by taking  $p$-th roots, namely,
\begin{equation}\label{eq:definition_Wp}
    \w_p(P,Q) = \min_{\pi \in \Pi(P,Q)} \left(\int_{\mathbb R^d \times \mathbb R^d} \|x-y\|^p d\pi(x,y) \right)^{1/p} \ .
\end{equation}
$\w_p$ is a metric on $\mathcal P_p(\mathbb R^d)$, the subset of $\mathcal P(\mathbb R^d)$ with finite $p$-th moment,  which metrizes the weak convergence of measures together with convergence of $p$-th moments. For these and other foundational results on optimal transport, we refer to \cite{villani2003Topics,villani2008optimal}.

Despite its intuitive definition and excellent properties, the Wasserstein distance presents a number of challenges for large-scale problems.  The first challenge is the computational complexity. 
Let $P_n$ and $Q_m$ be empirical measures from $n$ and $m$ i.i.d.\ observations from $P,\ Q\in\mathcal P_p(\mathbb R^d)$. 
If $n=m$, computing the Wasserstein distance between the empirical measures entails a cost that scales at a rate of $n\log(n)$ for $d=1$, using its closed-form expression as the $L^p$-norm between the quantile functions. For general $d>1$, no such simplification exists, and computing the Wasserstein distance requires solving a linear program,  which scales as $O(n^3)$ \cite{COTFNT}, limiting its applicability to domains with a heavy reliance on large-scale data. 
The second challenge arises from the significant degradation of the statistical performance of Wasserstein-based  analysis in high-dimensional settings. Understanding the convergence rates of empirical estimates of the distance to the theoretical values is a fundamental task. Under suitable moment assumptions, \cite{Fournier2015rateConvergence,Weed2019sharpAsymptotic} proved that if $P=Q$ and $d>2p$, then $\mathbb{E}(\w_p^p(P_n,Q_n))$ converges to zero at the rate $n^{-p/d}$. Similarly, if $P\neq Q$ and $d\geq 5$, assuming that both measures are absolutely continuous and sub-Weibull, \cite{manole2024sharpConvergence} showed that $\mathbb{E}(\w_p^p(P_n,Q_n)- \w_p^p(P,Q))$ converges at the rate $n^{-\min(2,p)/d}$, demonstrating that the Wasserstein distance suffers from the curse of dimensionality. 

This unfavorable asymptotic behavior hinders the derivation of central limit theorems (CLTs) for the Wasserstein distance in arbitrary dimensions, both for $P=Q$ and $P\neq Q$, when $P$ and $Q$ are sufficiently regular. The only available result in arbitrary dimension was first established in \cite{delBarrio2019generalDimension} for $p=2$, and later extended in \cite{delBarrio2024CLTgeneral} to a broader class of cost functions, including $c(x,y)=\|x-y\|^p$ for $p>1$,  yielding the CLT
\begin{equation}\label{eq:tcl_intro}
    \sqrt{n}\Bigl( \w_p^p(P_n,Q) - \mathbb{E}(\w_p^p(P_n,Q)) \Bigr)  \rightsquigarrow N(0,\sigma_{P,Q}^2) \ ,
\end{equation}  
where the Gaussian limit is non-degenerate if and only if $P\neq Q$. The convergence in \eqref{eq:tcl_intro} has limited applicability for inferential purposes, since the centering constant is not $\w_p^p(P,Q)$, as one would desire. Moreover, this convergence, along with the previous expectation bounds, ensures that a central limit theorem centered at the population value $\w_p^p(P,Q)$ cannot be established in arbitrary dimension if $P$ and $Q$ are sufficiently regular.
Consequently, distributional limits are mainly available in the discrete setting \cite{sommerfeld2017inferenceFiniteSpaces,tameling2019empirical} or in one dimension. In the latter, initial work focused on the case  $P=Q$ \cite{delBarrio1999,delBarrio2005asyptotics}, while for $P\neq Q$, \cite{delBarrio2019CLTfairness} provided sufficient conditions to ensure that the bias induced by replacing the expected value in \eqref{eq:tcl_intro} with $\w_p^p(P,Q)$ vanishes as $n\rightarrow \infty$, for $p\geq 2$. Alternative assumptions have also been explored; for example, \cite{hundrieser2024unifying} establishes a CLT for compactly supported measures in $d\le 3$.

To address these limitations, several alternative metrics have been proposed to preserve the favorable properties of the Wasserstein distance while alleviating its computational and statistical challenges. Notably, entropically regularized optimal transport  \cite{Cuturi2013SinkhornDL,COTFNT} adds a convex regularization term to the original problem, resulting in improved computational and statistical properties \cite{mena2019statisticalBounds,delBarrio2023improvedEntropic}. Numerous approaches have also been developed to leverage the unique aspects of the one-dimensional setting. In this work, we  focus in one of the most prominent proposals, the sliced Wasserstein distance \cite{Bonneel2015slicedAndRadon}, which is defined as an average over one-dimensional projections
\begin{equation}\label{defn:sliced}
\sw_p(P,Q) = \Bigl( \int_{\mathbb S^{d-1}} \w_p^p(P^\theta,Q^\theta) \ d\sigma(\theta) \Bigr)^{1/p},
\end{equation}
where $\operatorname{Pr}_\theta(x)=\langle \theta,x\rangle$ denotes the projection map, $P^\theta = \operatorname{Pr}_\theta \# P$, $Q^\theta = \operatorname{Pr}_\theta \# Q$ are the projected measures, and $\mathbb S^{d-1}$ is the unit sphere in $\mathbb R^d$. This distance preserves key properties of the Wasserstein distance while alleviating the curse of dimensionality. We now review related contributions in this direction.

\subsection*{Related work}

The statistical analysis of the sliced Wasserstein distance has evolved alongside other projection-based distances, such as the max-sliced distance, which replaces the integral in \eqref{defn:sliced} with the maximum over $\mathbb S^{d-1}$. The first important results for these projection-based distances were developed specifically for each problem. \cite{weed2022estimation} derived expectation bounds with rates independent of the dimension for the max-sliced version (and a subspace generalization). \cite{Manole2022minimax} adapted the previous bounds to the sliced setting, extending their validity to a trimmed version of \eqref{defn:sliced}, and proving the first distributional limit for the sliced Wasserstein distance and its trimmed version, for $p>1$ and $P\neq Q$ (in the sense that the asymptotic distribution vanishes if $P=Q$, as in \eqref{eq:tcl_intro}). Denoting by $F_{\theta,n}$ and $F_\theta$ the cumulative distribution function of the projected distributions, 
their asymptotic result is based on the weak convergence of the empirical process 
$\mathbb G_n(\theta,x) = \sqrt{n}(F_{\theta,n}(x)-F_{\theta}(x))$, along with Hadamard differentiability and the functional delta method. The work primarily focuses on the trimmed version and the assumptions are rather strong for the untrimmed setting, involving 
\begin{equation}\label{eq:SJ_infty}
    \operatorname{SJ}_{\infty}(P) = \underset{0<t<1}{\operatorname{essup}} \   \frac{1}{f_\theta(F_\theta^{-1}(t))} <\infty\ \ ,
\end{equation}
where $f_\theta$ denotes the density associated with $F_\theta$ and $\emph{essup}$ the essential supremum.
An important step in the asymptotic theory for projection-based distances was the unifying approach presented in \cite{xi2022distributional}. Using ideas similar to those in \cite{Manole2022minimax}, and inspired by the duality-based approach in \cite{hundrieser2024unifying}, \cite{xi2022distributional} investigates the weak convergence of the sliced process $\mathbb{G}_n(\theta) = \sqrt{n}(\w_p^p(P_n^\theta, Q^\theta) - \w_p^p(P^\theta, Q^\theta))$ for compactly supported probability measures such that the support of the projected probabilities is an interval. Weak convergence for various projection-based distance notions, including the sliced Wasserstein distance, follows directly from Hadamard differentiability and the functional delta method. Their main results primarily address  the case $p> 1$, $P \neq Q$, although some insights are also provided for $p=1$.

Independently, \cite{Goldfel2024statisticalInference} proved weak convergence of the empirical process indexed by functions $\varphi \circ \operatorname{Pr}_\theta$, where $\varphi$ is 1-Lipschitz and $\theta \in \mathbb{S}^{d-1}$. Using the dual formulation of $\w_1$ in one-dimension together with the extended functional delta method, distributional limits are obtained for $p>1$ and $P \neq Q$, under the assumption that the probability measures are compactly supported with convex support. Additionally, for $p=1$, distributional limits under mild assumptions are derived for both $P \neq Q$ and $P = Q$, with the latter assumptions subsequently weakened by \cite{xu2022central}.

Recently, \cite{hundrieser2024estimatedCosts} provided a refined version of the results in \cite{xi2022distributional} and \cite{Goldfel2024statisticalInference}, again for compactly supported probabilities but with slightly weaker assumptions on the supports, using similar arguments to those in \cite{xi2022distributional}. In a different vein, employing slightly different proof techniques, their work also explores, with great generality, distributional limits for empirical transport problems where the cost function is also estimated from the data.

\subsection*{Contributions}

This work establishes a new CLT for the empirical sliced Wasserstein distance, for $p>1$, based on the Efron-Stein inequality. Our methodology builds on \cite{delBarrio2019generalDimension}, later refined for one-dimensional distributions in \cite{delBarrio2019CLTfairness} and for general costs in \cite{delBarrio2024CLTgeneral}. The same ideas allowed \cite{mena2019statisticalBounds} and \cite{delBarrio2023improvedEntropic} to establish a CLT for the entropic transportation cost. Consequently, a natural question is whether these techniques can be adapted to the sliced Wasserstein distance; we provide an affirmative answer in Section \ref{sect:CLT_general}. The adaptation is not straightforward, with the main difficulty lying in the simultaneous treatment of optimal transport potentials across all directions $\theta \in \mathbb S^{d-1}$. To address this, we define a specific class of random optimal transport potentials exhibiting convergence and uniform integrability properties (see Section \ref{sect:OT_potentials}) which may be of independent interest. The precise control of these potentials underlies Theorem \ref{teor:CLT}, which ensures that
\begin{equation}\label{eq:tcl}
    \sqrt{n}\Bigl( \sw_p^p(P_n,Q) - \mathbb E(\sw_p^p(P_n,Q)) \Bigr)  \rightsquigarrow N(0,v^2_{P,Q}) \ 
\end{equation} 
under mild assumptions, where $v_{P,Q}^2 < \infty$ is a constant that vanishes when $P = Q$. A key contribution of Theorem \ref{teor:CLT} is that it enables the decomposition of the original problem into two distinct components:
(i) the analysis of the CLT centered at the expected empirical distance, as described in \eqref{eq:tcl}, and
(ii) the investigation of the bias term
\begin{equation}\label{eq:bias}
    \sqrt{n}\Bigl(\mathbb E(\sw_p^p(P_n,Q))  -  \sw_p^p(P,Q) \Bigr)  \ .
\end{equation} 
Section \ref{sect:bias} is devoted to a detailed study of this subject. In Theorem \ref{teor:bias_improved_csorgo_sliced}, we provide sufficient conditions ensuring that \eqref{eq:bias} converges to zero for general $p>1$, substantially simplifying existing one-dimensional results in \cite{delBarrio2019CLTfairness}. By combining these findings with Theorem \ref{teor:CLT}, we derive a CLT for the sliced Wasserstein distance centered at the population value, extending the asymptotic framework of \cite{Goldfel2024statisticalInference, hundrieser2024estimatedCosts,Manole2022minimax,xi2022distributional}.
A key novelty of our work is that, for the first time, the CLT applies to non-compactly supported measures, including the fundamental case of Gaussian distributions, central in theory and applications, representing a significant step forward in the statistical understanding of sliced Wasserstein distances beyond the compact-support setting.

Once the CLT centered at the population value is established, a further challenge for inference is the intractability of the integral defining $\sw_p^p(P_n, Q)$, which is typically approximated via Monte Carlo averaging over $k$ random directions. Section \ref{sect:slicing} addresses this issue: Proposition \ref{teor:finite_slicing} characterizes the asymptotic error, and Corollary \ref{coro_slicing} provides a CLT with a flexible asymptotic variance adapting to the limiting projection regime $\tau = \lim \frac{k}{k+n} \in [0,1]$. To our knowledge, this is the first CLT that explicitly captures the variability induced by the Monte Carlo approximation. Combined with the consistent estimation of the asymptotic variance (Section \ref{sect:variance}), this approach eliminates the need for bootstrap methods, avoiding the substantial computational overhead they typically entail. Finally, Section \ref{sect:simulations} presents a brief simulation study illustrating the practical impact of $\tau$. Full proofs and auxiliary results are provided in the appendix.

\subsection*{Notation}

Throughout this work, $\langle\cdot,\cdot\rangle$ and $\|\cdot\|$ denote the usual inner product and the Euclidean norm om $\mathbb R^d$. Generally, we  denote vectors in $\mathbb{R}^d$ by $x,y$, while reserving the notation $s,t$ for points in $\mathbb{R}$. $\mathbb S^{d-1}$ denotes the unit sphere, i.e., the set of vectors $\theta\in\mathbb R^d$ with $\|\theta\|=1$. We denote by $\mathcal P(\mathbb R^d)$ the set of probability measures on $\mathbb R^{d}$, and by $\mathcal P_p(\mathbb R^d)$ the subset of probabilities with finite $p$-th moment. Given $P \in \mathcal P(\mathbb R^d)$, $\supp{P}$ denotes the smallest closed subset of $\mathbb R^d$ such that $P(\supp{P})=1$. We say that $P$ has negligible boundary if 
\begin{equation*}
\ell_d\bigr(\supp{P}\setminus\operatorname{int}(\supp{P})\bigl) = 0 \ ,
\end{equation*}
where $\ell_d$ denotes the Lebesgue measure in $\mathbb R^d$. For each $\theta \in\mathbb S^{d-1}$, $P^\theta = \operatorname{Pr}_\theta \# P$ is the projected distribution onto the direction of $\theta$, where $\#$ is the push-forward operator and $\operatorname{Pr}_\theta(x)= \langle \theta,x\rangle$. The distribution function of $P^\theta$ is denoted by $F_\theta$, and the density, if it exists, by $f_\theta$. $X_1, \ldots, X_n$ are i.i.d. (independent and identically distributed) random variables drawn from $P$, and $Y_1, \ldots, Y_m$ i.i.d. random variables drawn from $Q$, assumed to be independent of the $X_i$'s. The uniform distribution on $\mathbb S^{d-1}$ is denoted by $\sigma$, and $\Theta_1,\ldots,\Theta_k$  represent i.i.d. random variables drawn from $\sigma$, independent of the $X_i$'s and $Y_i$'s. The probability measure induced by a random variable $Z$ is denoted by $\mathcal L(Z)$. The symbols $\rightarrow_p$ and $\rightsquigarrow$ are used to denote convergence in probability and weak convergence, respectively.  
We assume that all random variables are defined on a common probability space $(\Omega, \mathcal{F}, \mathbb{P})$, and given a random variable $X$ defined on $\Omega$, we denote the image of an element $\omega\in\Omega$ by $X^\omega$. Finally, additional section-specific notation is introduced as needed; in particular, $c$-concavity notation is presented in \Cref{subsec:prelim}.

\section{Sliced optimal transport maps and potentials}\label{sect:OT_potentials}

This section explains the duality theory that underlies the derivation of our main results. In the first part, we introduce the dual problem and review existing results for optimal transport maps and potentials. The second part focuses on adapting this theory to the sliced setting by considering a particular, carefully chosen set of optimal transport potentials across the values of $\theta\in \mathbb S^{d-1}$.

\subsection{Preliminaries}\label{subsec:prelim}

Despite the validity of the duality theory presented here for general cost functions, in this work, our attention will be restricted to the cost $c_p(x, y) = h_p(x - y) = \|x - y\|^p$ for $p > 1$. To simplify the notation, we usually denote the cost by $c = c_p$. With this convention, given $P, Q \in \mathcal{P}_p(\mathbb{R}^d)$, the Kantorovich duality guarantees that
\begin{equation}\label{eq:dual}
    \w_p^p(P,Q) = \max_{\phi\in L^1(P)} \int_{\mathbb R^d} \phi(x)dP(x) + \int_{\mathbb R^d} \phi^{c}(y)dQ(y) \ ,
\end{equation}
where $\phi^c$ is the $c$-conjugate of $\phi$, defined as $\phi^{c}(y) = \inf_{x\in\mathbb R^d} \bigl(\| x-y\|^p - \phi(x)\bigr)$. A maximizer in \eqref{eq:dual} is called a  \textit{optimal transport potential from $P$ to $Q$}. The theory of $c$-concavity plays a fundamental role, for which a thorough exposition can be found in \cite{gangbo1996geometry} and \cite{villani2003Topics}. In particular, we say that $\phi:\mathbb R^d \rightarrow [-\infty,\infty)$ is $c$-concave if there exists $\psi:\mathbb R^d \rightarrow [-\infty,\infty)$ such that $\phi=\psi^c$. Obviously, $\phi^c$ is $c$-concave, and it is easy to see that $\phi^{cc}\geq \phi$, where the inequality is an equality if and only if $\phi$ is $c$-concave. This implies that we can restrict the supremum in \eqref{eq:dual} to the subset of $c$-concave functions. Given a $c$-concave function $\phi:\mathbb R^d \rightarrow[-\infty,\infty)$, the $c$-superdifferential of $\phi$ is 
\begin{equation}\label{eq:defn_superdifferential}
    \partial^c \phi = \Big\{ (x,y)\in \mathbb R^d\times \mathbb R^d : f(z)\leq f(x) + \bigl(\|z-y\|^p-\|x-y\|^p\bigr) \ \forall \ z\in \mathbb R^d \Big\} \ . 
\end{equation}
For any $x \in \mathbb R^d$, we denote $\partial^c \phi(x) = \{y:(x,y)\in \partial^c \phi\}$. If $\partial^c \phi(x)$ consists of a single point, we denote $\nabla^c \phi(x) = \partial^c \phi(x)$. If $\phi$ is $c$-concave, $\phi(x)+\phi^c(y)\leq c(x,y)$, and the inequality is an equality if and only if $y\in \partial^c\phi(x)$. 
Theorems 2.3 and 2.7 in \cite{gangbo1996geometry} ensure that $\pi \in \Pi(P,Q)$ is optimal if and only if $\supp \pi \subset \partial^c \phi(x)$ for some $c$-concave function $\phi$, and in this case, $\phi$ is an optimal transport potential. Similarly, given any $c$-concave optimal transport potential $\phi$, there exists an optimal plan $\pi$ verifying $\supp \pi \subset \partial^c \phi(x)$. Moreover, Proposition 3.4 in \cite{gangbo1996geometry} guarantees that $\phi$ is differentiable at almost every point of $\operatorname{dom}(\phi)$, and if $P$ is absolutely continuous,  Theorem 1.2 in \cite{gangbo1996geometry} ensures that $\nabla^c \phi$ is the $P$-a.s. unique optimal transport map. In addition, if $P$ has negligible boundary and $\operatorname{int}(\supp P)$ is connected, Corollary 2.7 in \cite{delBarrio2024CLTgeneral} ensures the uniqueness, up to constants, of the optimal  potentials, and Theorem 3.4 in \cite{delBarrio2024CLTgeneral} establishes uniform convergence of empirical potentials, up to constants, on the compact sets of $\operatorname{int}(\supp P)$.

\subsection{Adaptation to the sliced setting}\label{sect:adaptation_potentials}

Let $P,Q\in\mathcal P_p(\mathbb R^d)$ be such that $P$ is absolutely continuous, with negligible boundary and $\operatorname{int}(\supp{P})$ connected, and let $P_n$ be the empirical distribution associated with $P$. Then, for every $\theta\in\mathbb S^{d-1}$, $P^\theta,Q^\theta \in \mathcal P_{p}(\mathbb R)$, and according to Lemma \ref{lemma:support}, $P^\theta$ is absolutely continuous with negligible boundary and $\operatorname{int}(\supp{P^\theta})$ connected. Therefore, in the set of $\mathbb P$-probability one satisfying $P_n\rightsquigarrow P$, Theorem 3.4 in \cite{delBarrio2024CLTgeneral} applies for every $\theta\in\mathbb S^{d-1}$. If we denote by $\phi^\theta$ and $\phi^\theta_n$ any $c$-concave optimal transport potential from $P^\theta$ to $Q^\theta$, and from $P^\theta_n$ to $Q^\theta$, respectively, there exist constants $a_n^\theta\in \mathbb R$ such that 
$\phi_n^\theta - a_n^\theta \rightarrow \phi^\theta$
in the sense of uniform convergence on the compact sets of $\operatorname{int}(\supp{P^\theta})$. To exploit duality within the sliced framework, we need a suitable choice of the constants
across the values of $\theta\in\mathbb S^{d-1}$. 

Inspired by \cite{Goldfel2024statisticalInference}, if we are only interested in the limit potentials $\phi^\theta$, we can proceed as follows. By the comments after Lemma \ref{lemma:c-convave},  given  $x_0\in\operatorname{int}(\supp P)$, it follows that $\phi^\theta(\langle \theta,x_0\rangle )\in\mathbb R$.  Thus, we can define 
\begin{equation}\label{defn:fixed_potentials_x0}
  \phi^\theta_{x_0} := \phi^\theta - \phi^\theta(\langle \theta,x_0\rangle ) \ .
\end{equation}
By Corollary 2.7 in \cite{delBarrio2024CLTgeneral},  $\phi_{x_0}^\theta$ is unique, in the sense that if $\Tilde \phi^\theta_{x_0}$ is any other $c$-concave optimal transport potential from $P^\theta$ to $Q^\theta$ such that $\Tilde \phi^\theta_{x_0}(\langle \theta,x_0\rangle)=0$, then $\phi_{x_0}^\theta(s) = \Tilde \phi^\theta_{x_0}(s)$ for every  $s\in\operatorname{int}(\text{supp} (P^\theta))$. 
However, this definition cannot be extended to empirical potentials $\phi_n^\theta$. If $\supp{Q^\theta}$ is not lower bounded (resp. upper bounded), and 
\begin{equation*}
    \langle\theta,x_0\rangle < \min_{i=1,\ldots,n} \langle\theta,X_i\rangle \quad \bigl(\textnormal{resp. }\langle\theta,x_0\rangle > \max_{i=1,\ldots,n} \langle\theta,X_i\rangle\big)
\end{equation*} 
then every $c$-concave optimal transport potential $\phi_n^\theta$ satisfies $\phi_n^\theta(\langle\theta,x_0\rangle)=-\infty$. To check this, consider an optimal transport plan $\pi_n^\theta$ such that  $\supp{\pi_n^\theta}\subset \partial^c \phi_n^\theta$. If $\supp{Q^\theta}$ is not lower bounded, by optimality of $\pi_n^\theta$, there exists $\{y_l\}_{l=1}^\infty$ such that $y_l\rightarrow -\infty$ and $( \min_{i=1,\ldots,n} \langle\theta,X_i\rangle ,y_l) \in \supp{\pi_n^\theta}\subset \partial^c \phi_n$. By  definition of $\partial^c \phi_n$, for every  $s< \min_{i=1,\ldots,n} \langle\theta,X_i\rangle$,
\begin{equation*}
    \phi_n^\theta(s) \leq \phi_n^\theta\Big(\min_{i=1,\ldots,n} \langle\theta,X_i\rangle\Big) + \Bigl( \big|s-y_l\big|^p - \big|\min_{i=1,\ldots,n} \langle\theta,X_i\rangle-y_l\big|^p\Bigr) \xrightarrow{l\rightarrow \infty} -\infty \ .
\end{equation*}
For $n$ large enough ($n$ depends on $\{X_n\}_{n=1}^\infty$), we can ensure $x_0\in [\min_{i=1,\ldots,n} \langle\theta,X_i\rangle,\max_{i=1,\ldots,n} \langle\theta,X_i\rangle]$, so $\phi^\theta_{n,x_0} = \phi_n^\theta - \phi_n^\theta(\langle \theta,x_0\rangle)$ is well defined. Elaborating on these ideas, one could prove convergence  $\phi_{n,x_0}^\theta \rightarrow \phi_{x_0}^\theta$. Nevertheless, this approach is not valid to provide integrability guarantees, which are crucial in our analysis.  
In contrast, we propose a simple  alternative solution to the problem. Since $X_1$ is the only point such that its projection belongs to  $[\min_{i=1,\ldots,n} \langle\theta,X_i\rangle,$ $\max_{i=1,\ldots,n} \langle\theta,X_i\rangle]$ for every $\theta\in\mathbb S^{d-1}$ and $ n\in\mathbb N$, we propose to consider $X_1$ as a reference point. This choice is justified by the following proposition.

\begin{proposition}\label{prop:regularidad_potenciales}
    Let $p>1$,  $c=c_p$, and $P,Q\in\mathcal P_p(\mathbb R^d)$ such that $P$ is absolutely continuous with negligible boundary and $\operatorname{int}(\supp{P})$ connected. For each $\theta\in\mathbb S^{d-1}$, let $\phi^\theta$, $\phi^\theta_n$ and  $\phi^\theta_{n,m}$ be any $c$-concave optimal transport potentials from $P^\theta$ to $Q^\theta$, from $P_n^\theta$ to $Q^\theta$ and from $P_n^\theta$ to $Q_m^\theta$, respectively.   
    Then, the following properties hold up to a set of $\mathbb P$-probability zero.  

       \begin{enumerate}[label=(\roman*)]
    \item For every $\theta\in\mathbb S^{d-1}$, 
    $
        \phi_{X_1}^\theta := \phi^\theta- \phi^\theta(\langle\theta,X_1\rangle)
    $
    is well defined, and it is a $c$-concave optimal transport potential from  $P^\theta$ to $Q^\theta$ satisfying $\phi_{X_1}^\theta(\langle\theta,X_1\rangle)=0$. Moreover, it is unique in  $\operatorname{int}(\supp{P^\theta})$.
    \item The map $\theta \mapsto \phi_{X_1}^\theta\left(\langle\theta, x\rangle\right)$ is continuous for every $x\in\operatorname{int}(\supp{P})$.

       \item For every $\theta\in\mathbb S^{d-1}$, $          \phi_{n,X_1}^\theta := \phi_n^\theta- \phi_n^\theta(\langle\theta,X_1\rangle) $ is well defined, and it is a $c$-concave optimal transport potential from  $P_n^\theta$ to $Q^\theta$ satisfying $\phi_{n,X_1}^\theta(\langle\theta,X_1\rangle)=0$.

        \item For every $\theta\in\mathbb S^{d-1}$, $\phi_{n,X_1}^\theta\rightarrow \phi_{X_1}^\theta$ in the sense of uniform convergence on the compacts of $\operatorname{int}(\supp{P^\theta})$,  and for each compact $K\subset \operatorname{int}(\supp{P^\theta})\cap \operatorname{dom}(\nabla\phi_{X_1}^\theta)$, 
    $
        \sup_{s\in K} \sup_{t_n\in\partial^c\phi_{n,X_1}^\theta(s)} |t_n-\nabla\phi_{X_1}^\theta(s)| \rightarrow 0 \ .
    $
    \item  Given $x \in \operatorname{int}(\operatorname{supp}(P))$, there exist  $n_0$ and $M_{x}<\infty$  such that $|\phi_{n,X_1}^\theta( \langle \theta, x\rangle)| \leq M_{x}$ for every  $n\geq n_0,\ \theta \in \mathbb S^{d-1}$.

        \item Let $m=m(n)\rightarrow\infty$ as $n\rightarrow \infty$. For every $\theta\in\mathbb S^{d-1}$, $
        \phi_{n,m,X_1}^\theta := \phi_{n,m}^\theta- \phi_{n,m}^\theta(\langle\theta,X_1\rangle)
    $
    is well defined, and it is a $c$-concave optimal transport potential from  $P_n^\theta$ to $Q_m^\theta$ verifying $\phi_{n,m,X_1}^\theta(\langle\theta,X_1\rangle)=0$. Moreover, assertions \textit{(iv)} and
    \textit{(v)} are also verified by $\phi_{n,m,X_1}^\theta$.
    \end{enumerate}

\end{proposition}

\noindent
Given any family of optimal transport potentials $\{\phi^\theta\}_{\theta\in\mathbb S^{d-1}}$, and any probability measure $\mu\in\mathcal P(\mathbb R^d)$, denote
\begin{equation*}
     \| \phi^\theta \|_{L^{q}(\sigma\times \mu)}^q := \int_{\mathbb S^{d-1}\times \mathbb R^d} \big| \phi^\theta (\langle\theta,x\rangle) \big|^q d(\sigma\times \mu)(\theta,x) \ .
\end{equation*}
Our next result establishes expectation bounds on the $L^q(\sigma\times\mu)$-norm of the $c$-concave random potentials $\phi^\theta_{X_1}$, $\phi^\theta_{n,X_1}$, $\phi^\theta_{n,m,X_1}$ defined as in \Cref{prop:regularidad_potenciales}, when  $\mu=P$ or $\mu=P_n$, and similarly for their $c$-conjugates. 

\begin{proposition}\label{prop:integrability_potentials}
Let $q\geq 1$. Under the assumptions of Proposition \ref{prop:regularidad_potenciales}, if  $P,Q$ have finite moments of order $pq$, then 
\begin{align}
     \mathbb E\Bigl( \| \phi_{X_1}^\theta \|_{L^{q}(\sigma\times P)}^{q}\Bigr)& \vee \   \mathbb E\Bigl( \| (\phi_{X_1}^\theta )^c\|_{L^{q}(\sigma\times Q)}^{q}\Bigr) < \infty  \ , \label{eq:int_potentials_limit}\\
      \sup_{n} \ \mathbb E\Bigl( \| \phi_{n,X_1}^\theta \|_{L^{q}(\sigma\times P_n)}^{q} \Bigr) & \vee \  \sup_{n} \ \mathbb E\Bigl( \| (\phi_{n,X_1}^\theta)^c \|_{L^{q}(\sigma\times Q)}^{q} \Bigr) < \infty \ , \label{eq:int_potentials_one_sample}\\
       \sup_{n,m} \ \mathbb E\Bigl( \| \phi_{n,m,X_1}^\theta \|_{L^{q}(\sigma\times P_n)}^{q} \Bigr) &\vee  \sup_{n,m} \ \mathbb E\Bigl( \| (\phi_{n,m,X_1}^\theta)^c \|_{L^{q}(\sigma\times Q_m)}^{q} \Bigr) 
 < \infty  \ .\label{eq:int_potentials_two_sample}
    \end{align}
\end{proposition}
\noindent
As an immediate consequence of \eqref{eq:int_potentials_limit}, since $\mathcal L(X_1)= P$,  we deduce that if we define the map  $g(x) = \|\phi_x^\theta\|_{L^{q}(\sigma\times P)}^{q}$  or $g(x) = \|(\phi_x^\theta)^c\|_{L^{q}(\sigma\times Q)}^{q}$ then $\int g(x) dP(x)<\infty$, which allows us to conclude the following result. 
\begin{corollary} \label{corollary:integrability}
 Under the same assumptions as in Proposition \ref{prop:integrability_potentials}, for almost every $x_0\in\operatorname{int}(\supp{P})$,
\begin{equation}\label{eq:int_potentials_limit_x0}
      \| \phi_{x_0}^\theta \|_{L^{q}(\sigma\times P)}^{q} \vee\       \| (\phi_{x_0}^\theta)^c \|_{L^{q}(\sigma\times Q)}^{q}  <\infty \ .
 \end{equation}
\end{corollary}

\section{Central limit theorem with centering constant}\label{sect:CLT_general}

The aim of this section is to prove a central limit theorem for the sliced Wasserstein distance, centered at the population value, as stated in Theorem~\ref{teor:CLT}. The argument relies on Propositions~\ref{teor:primal} and~\ref{teor:var_dual}, both derived from the Efron-Stein inequality. Proposition~\ref{teor:primal} exploits the primal formulation, whereas Proposition~\ref{teor:var_dual} employs the dual formulation (Sections~\ref{sec:primal} and~\ref{sec:dual}). Theorem~\ref{teor:CLT} then follows by standard probabilistic arguments.

\subsection{Efron-Stein and primal representation}\label{sec:primal}

Let  $Z_n=f(X_1, \ldots, X_n)$, where  $X_1, \ldots, X_n$ are independent random variables. Given an independent copy $(X_1', \ldots, X_n')$ of $(X_1, \ldots, X_n)$, denote by $Z_n^i=f(X_1, \ldots, X_i', \ldots, X_n)$, for each $i=1,\ldots,n$. The Efron-Stein inequality ensures that
$
    \operatorname{Var}(Z_n) \leq \tfrac{1}{2} \sum_{i=1}^n \mathbb E\left(Z_n-Z_n^i\right)^2=\sum_{i=1}^n \mathbb E\left(Z_n-Z_n^i\right)_{+}^2 \ .
$
In particular, if $X_1, \ldots, X_n$ are i.i.d. and $f$ is a symmetric function of $X_1, \ldots, X_n$, the bound simplifies to
\begin{equation}\label{eq:efron-stein}
    \operatorname{Var}(Z_n) \leq n \mathbb E\left(Z_n-Z_n'\right)_{+}^2 \ ,
\end{equation}
where $Z_n'=f\left(X_1', X_2, \ldots, X_n\right)$. 
Our first result applies the Efron--Stein inequality to the symmetric variable $Z_n=\sw_p^p(P_n,Q)$, with $P_n=\frac{1}{n}\sum_{i=1}^n \delta_{X_i}$, comparing it to $Z_n'=\sw_p^p(P_n',Q)$, where $P_n'=\frac{1}{n}(\delta_{X_1'}+\sum_{i=2}^n \delta_{X_i})$. Proposition~\ref{teor:primal} builds on the primal definition of the Wasserstein distance~\eqref{eq:definition_Wp}, constructing from the optimal coupling $\pi_n'^\theta$ between $P_n'^\theta$ and $Q^\theta$ a near-optimal coupling for $P_n^\theta$ and $Q^\theta$, yielding bounds on $Z_n-Z_n'$. Our proof is a straightforward adaptation of the ideas in Proposition 3.1 of \cite{delBarrio2019generalDimension} and Lemma 4.1 of \cite{delBarrio2024CLTgeneral} to the sliced setting, but differs by using optimal plans instead of maps, thereby removing smoothness assumptions on~$Q$ without the approximation argument of Corollary 4.3 in \cite{delBarrio2024CLTgeneral}.

\begin{proposition}[\textbf{Variance bound}]\label{teor:primal}
Let $p> 1$ and $P,Q \in \mathcal P_{2p}
(\mathbb R^d)$. Given independent random variables $X,X'$ with probability law $P$ and $Y$ with probability law $Q$, define 
$
    K_p(P,Q) := p^2 \mathbb E (   \| X -X'\|^{2p})^{1/p}  \mathbb E (  \|X-Y\|^{2p} ) ^{(p-1)/p} \ .
$
Similarly, define the constant  $K_p(Q,P)$ . Then, 
\begin{enumerate}
    \item[(a)] $\operatorname{Var} (\sw_p^p(P_n,Q)) \leq \frac{K_p(P,Q)}{n} $. 
 \item[(b)]
$\operatorname{Var} (\sw_p^p(P_n,Q_m)) \leq \frac{K_p(P,Q)}{n} + \frac{K_p(Q,P)}{m} $.
\end{enumerate}
\end{proposition}

\subsection{Efron-Stein and dual representation}\label{sec:dual}
Following \cite{delBarrio2024CLTgeneral,delBarrio2019generalDimension}, we exploit the dual representation of the Wasserstein distance~\eqref{eq:dual} to obtain an improved variance bound. Unlike the primal case, the adaptation to the sliced setting introduces significant technical challenges, due to the the simultaneous handling of all optimal transport potentials across the values of $\theta \in \mathbb S^{d-1}$. This is resolved using the results of \Cref{sect:adaptation_potentials}. In particular, for each $\theta \in \mathbb S^{d-1}$, let $\phi^\theta_n$ be a $c$-concave optimal transport potential from $P_n^\theta$ to $Q^\theta$, and assume $\phi^\theta = \phi^\theta_{x_0}$ as in \eqref{defn:fixed_potentials_x0}. Defining $\phi^\theta_{X_1}, \phi^\theta_{n,X_1}$ as in \Cref{prop:regularidad_potenciales}, the optimality of $\psi^\theta_{n} =(\phi^\theta_{n,X_1})^c \in L^1(Q)$ yields
\begin{align}
\w_p^p(P_n^\theta,Q^\theta) - \w_p^p(P_n'^\theta,Q^\theta)
&\leq \int_{\mathbb R} (\psi^\theta_{n})^c(s)d(P_n^\theta- P_n'^\theta)(s) =  \int_{\mathbb R} \phi^\theta_{n,X_1}(s)d(P_n^\theta- P_n'^\theta)(s) \notag \\
&= \int_{\mathbb R} (\phi^\theta_{n,X_1}(s)-\phi^\theta_{X_1}(s))d(P_n^\theta- P_n'^\theta)(s) + \int_{\mathbb R} \phi^\theta_{x_0}(s)d(P_n^\theta- P_n'^\theta)(s) \ , \label{eq:linearization}
\end{align}
where the last equality follows from the definition of $\phi^\theta_{X_1}$.
Applying this bound uniformly over $\mathbb{S}^{d-1}$ allows us to control the difference $Z_n - Z_n'$. Proposition~\ref{teor:var_dual} refines Proposition~\ref{teor:primal} by applying the Efron--Stein inequality to a modified variable $R_n$, obtained by subtracting the integral over $\mathbb{S}^{d-1}$ of the last term in \eqref{eq:linearization}. Thanks to the results of \Cref{sect:adaptation_potentials}, the proof then proceeds analogously to Theorem~3.2 in \cite{delBarrio2019generalDimension} and Theorem~4.6 in \cite{delBarrio2024CLTgeneral}.

\begin{theorem}\label{teor:var_dual}
    Let $p>1$ and $c=c_p$. Assume that $P,Q \in \mathcal P_{2p(1+\gamma)}(\mathbb R^d)$, for some $\gamma>0$, and that $P$ is absolutely continuous with negligible boundary and $\operatorname{int}(\supp{P})$ connected. Let $x_0\in \operatorname{int}(\supp P)$ satisfy \eqref{eq:int_potentials_limit_x0}. 
 
    \begin{itemize}
        \item[(a)] If we define
    $
        R_n := \sw_p^p(P_n,Q) - \int_{\mathbb S^{d-1}} \left( \int_{\mathbb R^d} \phi^\theta_{x_0}(\langle\theta,x\rangle) dP_n(x) \right) d\sigma(\theta) ,
    $
    then, as $n\rightarrow\infty$, $n \textnormal{Var}(R_n) \rightarrow 0$. 
    
    \item[(b)] Let $m=m(n)\rightarrow\infty$ as $n\rightarrow\infty$. If $Q$ is also absolutely continuous with negligible boundary and $\operatorname{int}(\supp{Q})$ is connected,     then, as $n\rightarrow\infty$, $\frac{nm}{n+m} \textnormal{Var}(R_{n,m}) \rightarrow 0$, where
    \begin{align*}
       \hspace*{-.6cm} R_{n,m} := \sw_p^p(P_n,Q_m) &- \int_{\mathbb S^{d-1}} \left( \int_{\mathbb R^d} \phi_{x_0}^\theta(\langle\theta,x\rangle) dP_n(x) \right) d\sigma(\theta)-\int_{\mathbb S^{d-1}} \left( \int_{\mathbb R^d} (\phi_{x_0}^\theta)^c(\langle\theta,y\rangle) dQ_m(y) \right) d\sigma(\theta) \ .
    \end{align*}
    \end{itemize}
\end{theorem}

\subsection{Central limit theorem}\label{subsec:clt}

Hereafter, given any $c$-concave optimal transport potential $\phi^\theta$ and $x\in\mathbb R^d$, denote, with a slight abuse of notation,   $\phi^\theta(x) = \phi^\theta(\langle \theta, x\rangle)$.
Under the assumptions of Proposition \ref{teor:var_dual}, applying Chebyshev's inequality,
\begin{equation} \label{eq:decomposition}
  \sqrt{n}(R_n - \mathbb E(R_n))  = \sqrt{n} \Bigl(\sw_p^p(P_n,Q) - \mathbb E\bigl(\sw_p^p(P_n,Q)\bigr)\Bigr) -  \sqrt{n}\Bigl( \frac{1}{n} \sum_{i=1}^n Y_i - \mathbb E(Y_i)\Bigr) \rightarrow_p 0 \ ,
\end{equation}
where $Y_i =  \int_{\mathbb S^{d-1}}   \phi^\theta_{x_0}(X_i) d\sigma(\theta)$,  $i=1,\ldots,n$, are i.i.d. random variables  with finite variance, by \eqref{eq:int_potentials_limit_x0}. Thus, applying the CLT to these variables together with \eqref{eq:decomposition} yields:

\begin{theorem}\label{teor:CLT}
Let $p>1,\gamma>0$ and $P,Q \in \mathcal P_{2p(1+\gamma)}(\mathbb R^d)$ be such that $P$ is absolutely continuous with negligible boundary  and $\operatorname{int}(\supp{P})$  connected. For each $\theta\in\mathbb S^{d-1}$, let $\phi^\theta$ be any $c$-concave optimal transport potential from $P^\theta$ to $Q^\theta$, and define $v_{P,Q}^2=  \int_{\mathbb S^{d-1}} \int_{\mathbb S^{d-1}} \operatorname{Cov}_P(\phi^\theta,\phi^\eta) d\sigma(\theta) d\sigma(\eta)$, and similarly  $v_{Q,P}^2$. Then:

\begin{itemize}
    \item[(a)] 
    $
        \sqrt{n} \Bigl(\sw_p^p(P_n,Q) - \mathbb E\bigl(\sw_p^p(P_n,Q)\bigr)\Bigr) \rightsquigarrow N(0,v_{P,Q}^2) \ .
    $
    \item[(b)]   If $m=m(n)\rightarrow\infty$ as $n\rightarrow\infty$ is such that $\frac{n}{n+m}\rightarrow \lambda \in [0,1]$ and $Q$ is also absolutely continuous with negligible boundary and $\operatorname{int}(\supp{Q})$  connected, 
    \begin{equation*}
        \sqrt{\frac{nm}{n+m}} \Bigl(\sw_p^p(P_n,Q_m) - \mathbb E\bigl(\sw_p^p(P_n,Q_m)\bigr)\Bigr) \rightsquigarrow N\bigl(0,(1-\lambda)v_{P,Q}^2+\lambda v_{Q,P}^2\bigr) \ .
    \end{equation*}
\end{itemize}
Moreover, in both cases, convergence in the stronger $\w_2$ sense also holds.
\end{theorem}

\section{Changing the centering constant}\label{sect:bias}

In this section, we study the problem of replacing the expected empirical values in \Cref{teor:CLT} with the corresponding population quantities. \Cref{teor:bias_improved_csorgo_sliced} provides sufficient assumptions ensuring that  
\begin{align}\label{eq:bias_conv_0}
0 \leq \sqrt{n} \Bigl( \mathbb{E}\bigl(\sw_p^p(P_n, Q)\bigr) - \sw_p^p(P, Q) \Bigr) \rightarrow 0 \  .
\end{align}
The non-negativity follows from the same duality argument as in \cite{delBarrio2019CLTfairness}.  
The proof of \Cref{teor:bias_improved_csorgo_sliced} relies on the following lemma, which establishes  weak convergence of weighted integrals of the empirical quantile process. The result is closely related to Theorem~2.1 in \cite{csorgo1990distributions}, which, as noted on page~67, also holds for \(p = 1\) without the absolute value. Building on the ideas of \cite{delBarrio2005asyptotics}, our proof relaxes the integrability conditions required for \(q(t) = f(F^{-1}(t)) / h(t)\) (connecting the notation of \cite{csorgo1990distributions} and \Cref{lemma:csorgo}).  
For clarity, the assumptions are stated separately for \(f \circ F^{-1}\) and \(h\).

\begin{lemma}\label{lemma:csorgo}
    Let $P \in \mathcal P(\mathbb R)$, with density $f$ such that $f(F^{-1}(t))$ is positive and continuous in $(0,1)$, and monotone for $t$ sufficiently close to $0$ and $1$. Let $h:\mathbb (0,1) \rightarrow \mathbb R$ be such that there exists $\alpha,\beta>1$, $\frac{1}{\alpha}+\frac{1}{\beta}=1$, verifying  
    \begin{itemize}
        \item[(i)] $n^{-\alpha/2}\int_{1/n}^{1-1/n} \frac{1}{f(F^{-1}(t))^\alpha} dt \rightarrow 0$ , 
        \item[(ii)] $ \int_0^1\int_0^1  \bigl(\frac{s\wedge t-st}{f(F^{-1}(s))f(F^{-1}(t))}\bigr)^\alpha ds dt<\infty$ 
        \item[(iii)] $\int_0^1 |h(t)|^\beta dt<\infty$.
    \end{itemize}
    Then, $ \int_{1/n}^{1-1/n} h(t) \sqrt{n}( F_n^{-1}(t) - t ) dt \rightsquigarrow N(0,v^2)$, where $v^2 := \int_0^1\int_0^1  \frac{h(s)h(t)}{f(F^{-1}(s))f(F^{-1}(t))}(s\wedge t-st) dtds<\infty$. Moreover, convergence also holds replacing $\int_{1/n}^{1-1/n}$ by $\int_0^1$ if either (i) and  (ii) are replaced by the stronger assumption $J_\alpha(P):= \int_0^1 \frac{(t(1-t))^{\alpha/2}}{f(F^{-1}(t))^{\alpha}}  dt<\infty$, or if there exists  $p>1$ such that $P\in\mathcal P_{2p}(\mathbb R)$ and $\int_0^1 |h(t)|^{2p/(p-1)} dt<\infty$.
\end{lemma}

The proof of this lemma proceeds in two steps. First, following the approach of \cite{delBarrio2005asyptotics}, we show that convergence holds when the general quantile process, appropriately weighted, is replaced by the uniform quantile process in the integral. Then, combining the ideas of \cite{csorgo1990distributions} with the previous argument, we demonstrate that this approximation does not affect the asymptotic distribution. Consequently, \Cref{lemma:csorgo} yields weaker assumptions than those that would result from Theorem~2.1 in \cite{csorgo1990distributions} together with H\"older’s inequality.
Moreover, by incorporating \Cref{lemma:csorgo} into the proof of Theorem~2.3 in \cite{delBarrio2019CLTfairness}, one can obtain simpler and more natural assumptions to establish \eqref{eq:bias_conv_0} for $d = 1$ and $p \geq 2$, the setting considered there. In particular, to prove the convergence of the integral of the uniform quantile process, \cite{delBarrio2019CLTfairness} follows a similar strategy, relying on Assumptions~(i) and~(ii) in \Cref{lemma:csorgo}, which are slightly weaker than the condition $J_\alpha(P) < \infty$.
However, to justify the approximation, they impose strong smoothness requirements—specifically, twice differentiability of the distribution function $F$ of $P$, and
$\sup_{t \in (0,1)} t(1-t) \, {|f'(F^{-1}(t))|f^{-2}(F^{-1}(t))} < \infty$---
typically imposed to ensure the asymptotic equivalence as random processes in $L^2(0,1)$, which is considerably stronger than the integral-level equivalence required in our analysis. Instead, we assume monotonicity in the tails, which---although not directly comparable---yields a simpler and more natural condition specifically tailored for the integral-level equivalence.\\

Our extension to the sliced setting builds upon the one-dimensional case. However, to prove \eqref{eq:bias_conv_0}, additional control across projections is required. In particular, for each \(\theta \in \mathbb{S}^{d-1}\), we replace Assumptions~(i) and~(ii) in \Cref{lemma:csorgo} for \(P_\theta\) with the stronger condition \(\operatorname{J}_\alpha(P_\theta) < \infty\). Specifically, for \(P \in \mathcal{P}(\mathbb{R}^d)\), we assume finiteness of
$\operatorname{SJ}_\alpha(P) := \int_{\mathbb{S}^{d-1}} \operatorname{J}_\alpha(P^\theta) \, d\sigma(\theta)$, originally introduced by \cite{Manole2022minimax}.
This assumption plays a crucial role in the proof, as it allows us to adopt a different approach from that of \cite{delBarrio2019CLTfairness} in the one-dimensional case, based on a careful analysis of convexity in \Cref{lemma:convergence_derivative}. This argument enables us to consider any $p > 1$, whereas Proposition~2.3 in \cite{delBarrio2019CLTfairness} applies only to $p \geq 2$, and eliminates the need for their assumption~(2.7), while still avoiding the strong twice-differentiability assumptions  imposed in their approach.

\begin{theorem}\label{teor:bias_improved_csorgo_sliced}
Assume $p>1$ and let $P,Q\in \mathcal P_{2p}(\mathbb R^d)$ be such that for every $\theta\in\mathbb S^{d-1}$, $P^\theta$ has distribution function $F_\theta$ and density $f_\theta$ verifying that $f_\theta(F_\theta^{-1}(t))$ is positive and continuous in $(0,1)$, and monotone for $t$ sufficiently close to $0$ and $1$.  If there exist $\alpha,\beta>1 $ with $\frac{1}{\alpha}+\frac{1}{\beta}=1$  such that $\operatorname{SJ}_{\alpha}(P) < \infty$ and $Q$, if $p=2$, or $P,Q$, if $p\neq 2$, have finite moment of order $\beta\max(1,p-1)$, 
then \eqref{eq:bias_conv_0} holds. Similarly, if $m=m(n)\rightarrow \infty$ as $n\rightarrow \infty$ and there exist $\alpha',\beta'>1$ such that all the previous assumptions are also verified exchanging the roles of $P$ and $Q$, then, as $n\rightarrow \infty$,
\begin{align}\label{remove_bias_W2_two_sample}
    \sqrt{\frac{nm}{n+m}} \Bigl( \mathbb E\bigl(\sw_p^p(P_n,Q_m)\bigr) - \sw_p^p(P,Q)\Bigr) \rightarrow 0 \ .
\end{align}
\end{theorem}

\noindent
As a straightforward consequence of Theorems  \ref{teor:CLT} and \ref{teor:bias_improved_csorgo_sliced}, we can  derive a CLT centered at $\sw_p^p(P,Q)$. 

\begin{corollary}\label{coro:CLT}
    Let $p>1$, and $P,Q \in \mathbb R^d$. If the assumptions of Theorems  \ref{teor:CLT} and \ref{teor:bias_improved_csorgo_sliced} for one-sample hold, then 
    \begin{equation} \label{eq:CLT_sliced_one_sample}
              \sqrt{n} \Bigl( \mathbb \sw_p^p(P_n,Q) - \sw_p^p(P,Q)\Bigr) \rightsquigarrow N(0,v_{P,Q}^2) \ .  
    \end{equation}
    Moreover, if the assumptions of the two-sample setting are also verified, then 
\begin{equation}
\label{eq:CLT_sliced_two_sample}
     \sqrt{\frac{nm}{n+m}} \Bigl(\sw_p^p(P_n,Q_m) -  \sw_p^p(P,Q)\Bigr) \rightsquigarrow N\bigl(0,(1-\lambda)v_{P,Q}^2+\lambda v_{Q,P}^2\bigr)  \ .
\end{equation}
\end{corollary}

\Cref{coro:CLT} extends and refines existing results. It requires the existence of $\alpha > 1$ such that $\operatorname{SJ}_\alpha(P) < \infty$ and sufficiently many finite moments for both $P$ and $Q$. 
For instance, if $P$ is a multivariate Gaussian distribution, then all moments are finite, and $\operatorname{SJ}_\alpha(P) < \infty$ if and only if $\alpha < 2$. 
Due to the trade-off between $\alpha$ and $\beta$, it follows that \eqref{eq:CLT_sliced_one_sample} and \eqref{eq:CLT_sliced_two_sample} hold for all $p > 1$ when $P,Q$ are multivariate Gaussian distributions, which was only known for the one-dimensional case and $p \geq 2$  \cite{delBarrio2019CLTfairness}. 
This also highlights the novelty of \Cref{coro:CLT}, which provides the first central limit theorem for $\sw_p^p$, $p > 1$, valid for measures without compact support and encompassing the fundamental Gaussian case. 
Moreover, while \cite{Manole2022minimax} conjectured that a CLT of the form \eqref{eq:CLT_sliced_one_sample} would require $\operatorname{SJ}_p(P) < \infty$, Corollary~\ref{coro:CLT} shows that it suffices to assume $\operatorname{SJ}_\alpha(P) < \infty$ for some $\alpha > 1$, provided that $P$ and $Q$ have sufficiently many moments.

\section{Slicing with a finite number of directions}\label{sect:slicing}

Theorem \ref{teor:bias_improved_csorgo_sliced} establishes a CLT centered at the population value $\sw_p^p(P,Q)$. However, to carry out inference on this quantity, we must evaluate the integral over $\mathbb{S}^{d-1}$ in the definition of $\sw_p^p(P_n,Q)$. In practice, this integral is typically approximated using the Monte Carlo estimator based on an independent sample $\Theta_1,\ldots,\Theta_k$ drawn  from $\sigma$, given by
\begin{equation}\label{eq:defn_finite_sliced} \sw_{p,k}^p(P_n,Q) := \frac{1}{k}\sum_{i=1}^{k} \w_p^p(P_n^{\Theta_i},Q^{\Theta_i}) \ . \end{equation}
Thus, an additional source of variation is introduced into the problem, which is controlled by the following proposition. 
\begin{proposition}\label{teor:finite_slicing}
    Let $p>1$, and $P,Q \in \mathcal P_{2p}(\mathbb R^d)$. For each $k\in\mathbb N$, let $\Theta_1,\ldots,\Theta_{k}$ be i.i.d. samples drawn from $\sigma$. If $k=k(n),m=m(n)\rightarrow \infty$ as $n\rightarrow \infty$, then
    \begin{equation}
        \begin{aligned}       
        S_{n,k}=& \sqrt{k}\Bigl( \sw_{p,k}^p(P_n,Q)-\sw_p^p(P_n,Q)\Bigr) &\rightsquigarrow N(0,w_{P,Q}^2) \ ,   \\
        S_{n,m,k}=& \sqrt{k}\Bigl( \sw_{p,k}^p(P_n,Q_m)-\sw_p^p(P_n,Q_m)\Bigr) &\rightsquigarrow N(0,w_{P,Q}^2)  \ , \label{eq:CLT_slicing_two_sample}
    \end{aligned}
    \end{equation}
    where $w_{P,Q}^2= \int_{\mathbb S^{d-1}} \w_p^{2p}(P,Q)d\sigma(\theta) - \sw_p^{2p}(P,Q)$. In fact, convergence in the stronger $W_2$ sense also holds. 
\end{proposition}

\noindent
Combining Theorems \ref{teor:CLT} and \ref{teor:bias_improved_csorgo_sliced}, and Proposition \ref{teor:finite_slicing}, we obtain:

\begin{corollary}\label{coro_slicing}
    If $P,Q\in\mathcal P_{2p}(\mathbb R^d)$ verify all the assumptions of the one-sample setting in Theorems \ref{teor:CLT} and  \ref{teor:bias_improved_csorgo_sliced}, and $k=k(n)\rightarrow\infty$ as $n\rightarrow \infty$, such that $k/(k+n) \rightarrow \tau \in[0,1]$, then 
    \begin{align}
         \sqrt{\frac{kn}{k+n}} \Bigl( \sw_{p,k}^p&(P_n,Q)-\sw_p^p(P,Q)\Bigr) \rightsquigarrow N\Bigl(0,(1-\tau)w_{P,Q}^2+\tau v_{P,Q}^2\Bigr) \ .
         \label{TCL:finite_sliced_one_sample}
    \end{align}
    Moreover, if the assumptions of the two sample-setting are also verified and $ k  / (k+ {\frac{nm}{n+m}})\rightarrow \tau\in[0,1]$, then 
    \begin{align}
        \sqrt{\frac{k {\frac{nm}{n+m}}}{k+ {\frac{nm}{n+m}}}} \Bigl( \sw_{p,k}^p&(P_n,Q_m)-\sw_p^p(P,Q)\Bigr)  \rightsquigarrow N\Bigl(0,(1-\tau)w_{P,Q}^2 + \tau \bigl((1-\lambda)v_{P,Q}^2+\lambda \label{TCL:finite_sliced_two_sample}v_{Q,P}^2\bigr)\Bigr) \ .
    \end{align}
\end{corollary}

\section{Variance estimation}\label{sect:variance}

The purpose of this section is to provide consistent estimators of the asymptotic variances. The simplest problem is the estimation of $w_{P,Q}^2$, defined in \Cref{teor:finite_slicing}. With the previous notation:

\begin{proposition}
\label{coro:consistW}
    Define  $    \hat w_{P_n,Q}^2 := \frac{1}{k} \sum_{i=1}^k \w_p^{2p}(P_n^{\Theta_i},Q^{\Theta_i}) - \sw_{p,k}^{2p}(P_n,Q)$ and similarly $\hat w_{P_n,Q_m}^2$ with $Q$ replaced by $Q_m$. Then, under the assumptions of Proposition \ref{teor:finite_slicing}, $\hat w_{P_n,Q}^2\rightarrow_p w_{P,Q}^2$ and $\hat w_{P_n,Q_m}^2 \rightarrow_p w_{P,Q}^2 \ .$
\end{proposition}

Estimating $v_{P,Q}^2$ (or symmetrically $v_{Q,P}^2$) is more delicate. Our approach builds on the expression in Theorem \ref{teor:CLT}, replacing $P$ and $\phi^\theta$ with their empirical counterparts and approximating the integral via Monte Carlo. The proof of \Cref{prop:consistencia} again relies on the control of the optimal transport potentials across directions $\theta \in \mathbb{S}^{d-1}$, developed in \Cref{sect:adaptation_potentials}, together with the strong law of large numbers for generalized U-statistics.

\begin{proposition}\label{prop:consistencia}
    Let $\phi_n^\theta$ (resp. $\phi_{n,m}^\theta$) denote any $c$-concave optimal transport potential from $P_n^\theta$ to $Q^\theta$ (resp. $Q_m^\theta$). Under the assumptions of Theorem \ref{teor:CLT} for the one-sample setting, $\hat v_{P_n,Q}^2 := \frac{1}{k^2}  \sum_{i,j=1}^{k} \operatorname{Cov}_{P_n}(\phi_n^{\Theta_i},\phi_n^{\Theta_j}) \rightarrow_p v_{P,Q}^2$. Moreover, under the assumptions of the two-sample setting, $\hat v_{P_n,Q_m}^2 := \frac{1}{k^2}  \sum_{i,j=1}^{k} \operatorname{Cov}_{P_n}(\phi_{n,m}^{\Theta_i},\phi_{n,m}^{\Theta_j})\rightarrow_p v_{P,Q}^2$. 
\end{proposition}

\section{Simulations}\label{sect:simulations}

\begin{figure}[h!]
    \centering
\includegraphics[width=.8\linewidth]{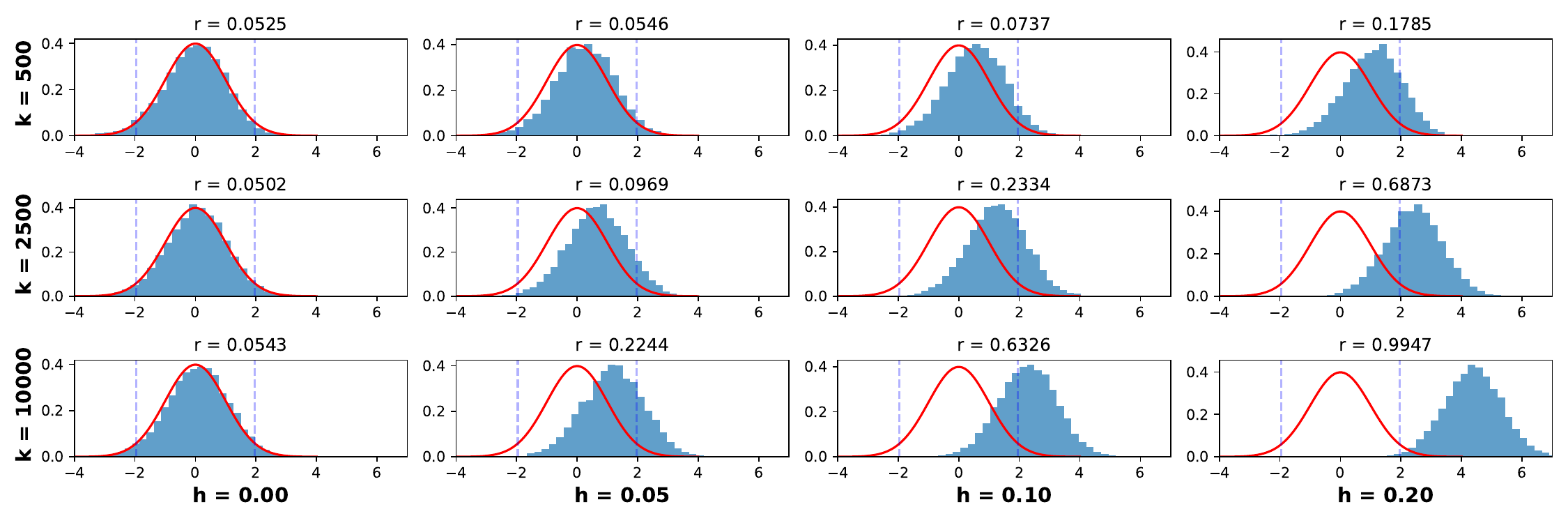}
    \caption{Histograms based on $10000$ independent samples of \eqref{CLT:inference} for $p=2$, $d=32$, $n=5000$, $m=3000$, $P=N(0,I_d)$, and $Q=Q_h=N((\sqrt{d}+h,0,\ldots,0),I_d)$, with $\Delta = \sw_2^2(P,Q_0)=1$, for different values of $k$ and $h$. In each subplot, $r$ denotes the proportion of samples satisfying $|T_{n,m,k}^\Delta| > q_{0.975}$, where $q_{0.975}$ is the $0.975$ quantile of a standard normal $N(0,1)$. Vertical lines represent $\pm q_{0.975}$.}
    \label{fig:histograms}
\end{figure}

This section demonstrates the applicability of our method for statistical inference and examines the role of the random projection regime, determined by $\tau = \lim k / (k+ {\tfrac{nm}{n+m}})\in[0,1]$. We focus on the two-sample setting because of its prominent importance in statistical applications. Provided the assumptions stated in the previous sections hold, one can define a consistent estimator $\hat v^2$ of the asymptotic variance in \eqref{TCL:finite_sliced_two_sample}, and given $\Delta>0$, if we define 
\begin{align}\label{CLT:inference}
 T_{n,m,k}^\Delta :=\frac{1}{\hat v}\sqrt{\frac{k\frac{nm}{n+m}}{k+\frac{nm}{n+m}}} \Bigl( \sw_{p,k}^p&(P_n,Q_m)-\Delta\Bigr) \ ,
\end{align}
then  $T_{n,m,k}^{\Delta} \rightsquigarrow N(0,1)$ under $H_0:\sw_p^p(P,Q) = \Delta$.
This provides an effective tool to perform  inference on the value $\sw_p^p(P,Q)$, which is valid for any value of $\tau\in[0,1]$. However, for inferential purposes, the regime with high number of projections ($\tau=1$),  is particularly relevant. 
Taking a closer look at the proof of Corollary \ref{coro_slicing}, we observe that, under the null hypothesis,  $\hat v T_{n,m,k}^\Delta$ behaves in the limit as $\sqrt\tau \times \eqref{eq:CLT_sliced_two_sample} + \sqrt{1-\tau} \times  \eqref{eq:CLT_slicing_two_sample}$. However, under $H_a:\sw_p^p(P,Q) \neq  \Delta$, $\hat v T_{n,m,k}^\Delta$ behaves in the limit as $\sqrt\tau \times \sqrt{\frac{nm}{n+m}} \bigl(\sw_p^p(P_n,Q_m) -  \Delta\bigr) + \sqrt{1-\tau} \times  \eqref{eq:CLT_slicing_two_sample}$. While the first term only converges to a Gaussian distribution under $H_0$, the convergence in  \eqref{eq:CLT_slicing_two_sample} does not depend on the value of $\Delta$. Therefore,  to achieve higher statistical power under the alternative hypothesis, more weight should be assigned to the divergent part, i.e., $\tau$ should be close to $1$. Alternatively, large values of $k$ are required. This behavior is illustrated in Figure \ref{fig:histograms}, 
where using Monte Carlo simulations, the distribution of $T_{n,m,k}^\Delta$ has been approximated under the null hypothesis and different alternatives in a Gaussian example.  
We can conclude as expected. The statistical power of the test under $H_a$ decreases significantly for small values of $k$. Additional computational details can be found in Appendix \ref{appendix:efficient}.

\section*{Conclusion}

This work provides novel contributions to the statistical theory of optimal transport, in particular to the asymptotic understanding of the sliced Wasserstein distance \(\sw_p^p\) for \(p>1\). We establish, for the first time, a CLT valid for measures without compact support, encompassing the fundamental case of Gaussian distributions. Our approach leverages the Efron-Stein inequality and introduces a rigorous framework for handling random optimal transport potentials across all directions \(\theta \in \mathbb{S}^{d-1}\), which enables the derivation of a CLT centered at the expected empirical value.
Beyond extending one-dimensional results \cite{delBarrio2019CLTfairness} to the sliced setting, we simplify the assumptions required to remove the centering constant and prove their validity for all \(p>1\).
Additionally, we propose a novel approach that incorporates the Monte Carlo approximation error directly into the limit theorem, explicitly quantifying its asymptotic effect, and provide a consistent estimator of the asymptotic variance. Together, these results yield a complete and computationally efficient inference framework for \(\sw_p^p\) that avoids the need for bootstrap methods.

An important direction for future research is to identify necessary and sufficient conditions under which the bias in \eqref{eq:bias_conv_0} is \(O(1)\) for \(P \neq Q\), in the spirit of \cite{bobkov2019onedimensional}, thereby refining the bounds of \cite{Manole2022minimax}. Moreover, as is common in optimal transport, the degenerate case \(P = Q\) remains particularly challenging, since the limiting variance vanishes and alternative techniques are required. Developing a comprehensive inference framework that also accommodates this case constitutes an exciting open problem.


\section*{Acknowledgments}

The first and second authors were supported by MCIN/AEI/10.13039/501100011033/FEDER under grant  PID2021-128314NB-I00. The third author was supported by the Agence Nationale de la Recherche under grant ANR-23-CE23-0029 Regul-IA. He also acknowledges the support from the AI Cluster ANITI (ANR-19-PI3A-0004).

\bibliographystyle{alpha}
\bibliography{references}  

@article{delBarrio2019generalDimension,
author = {Eustasio del Barrio and Jean-Michel Loubes},
title = {{Central limit theorems for empirical transportation cost in general dimension}},
volume = {47},
journal = {The Annals of Probability},
number = {2},
publisher = {Institute of Mathematical Statistics},
pages = {926 -- 951},
year = {2019},
doi = {10.1214/18-AOP1275}
}

@article{delBarrio2019CLTfairness,
    author = {del Barrio, Eustasio and Gordaliza, Paula and Loubes, Jean-Michel},
    title = {{A central limit theorem for $L^p$ transportation cost on the real line with application to fairness assessment in machine learning}},
    journal = {Information and Inference: A Journal of the IMA},
    volume = {8},
    number = {4},
    pages = {817-849},
    year = {2019},
    issn = {2049-8772},
    doi = {10.1093/imaiai/iaz016}
}

@article{delBarrio2024CLTgeneral,
author = {Eustasio del Barrio and Alberto Gonz{\'a}lez-Sanz and Jean-Michel Loubes},
title = {{Central limit theorems for general transportation costs}},
volume = {60},
journal = {Annales de l'Institut Henri Poincaré, Probabilités et Statistiques},
number = {2},
publisher = {Institut Henri Poincaré},
pages = {847 -- 873},
year = {2024},
doi = {10.1214/22-AIHP1356}
}

@article{Goldfel2024statisticalInference,
    author = {Goldfeld, Ziv and Kato, Kengo and Rioux, Gabriel and Sadhu, Ritwik},
    title = {{Statistical inference with regularized optimal transport}},
    journal = {Information and Inference: A Journal of the IMA},
    volume = {13},
    number = {1},
    pages = {iaad056},
    year = {2024},
    doi = {10.1093/imaiai/iaad056}
}

@article{Manole2022minimax,
author = {Tudor Manole and Sivaraman Balakrishnan and Larry Wasserman},
title = {{Minimax confidence intervals for the Sliced Wasserstein distance}},
volume = {16},
journal = {Electronic Journal of Statistics},
number = {1},
publisher = {Institute of Mathematical Statistics and Bernoulli Society},
pages = {2252 -- 2345},
year = {2022},
doi = {10.1214/22-EJS2001}
}

@article{manole2024sharpConvergence,
author = {Tudor Manole and Jonathan Niles-Weed},
title = {{Sharp convergence rates for empirical optimal transport with smooth costs}},
volume = {34},
journal = {The Annals of Applied Probability},
number = {1B},
publisher = {Institute of Mathematical Statistics},
pages = {1108 -- 1135},
year = {2024},
doi = {10.1214/23-AAP1986}
}

@article{Fournier2015rateConvergence,
author = {Fournier, Nicolas and Guillin, Arnaud},
title = {{On the rate of convergence in Wasserstein distance of the empirical measure}},
journal = {Probability Theory and Related Fields},
year = {2015},
volume = {162},
pages = {707–738},
doi = {https://doi.org/10.1007/s00440-014-0583-7}
}

@article{hundrieser2024estimatedCosts,
title = {{Empirical optimal transport under estimated costs: Distributional limits and statistical applications}},
journal = {Stochastic Processes and their Applications},
volume = {178},
pages = {104462},
year = {2024},
issn = {0304-4149},
doi = {https://doi.org/10.1016/j.spa.2024.104462},
author = {Shayan Hundrieser and Gilles Mordant and Christoph A. Weitkamp and Axel Munk}
}

@inproceedings{xi2022distributional,
 author = {Xi, Jiaqi and Niles-Weed, Jonathan},
 booktitle = {Advances in Neural Information Processing Systems},
 pages = {13961--13973},
 publisher = {Curran Associates, Inc.},
 title = {{Distributional convergence of the sliced Wasserstein process}},
 volume = {35},
 year = {2022}
}

@article{bobkov2019onedimensional,
title = {{One-dimensional empirical measures, order statistics, and Kantorovich transport distances}},
journal = {Memoirs of the American Mathematical Society}, 
volume = {261},
number = {1259}, 
author = {Bobkov, Sergey G. and Ledoux, Michel},
year = {2019}
}

@article{csorgo1978strongApproximations,
author = {Miklos Csorgo and Pal Revesz},
title = {{Strong approximations of the quantile process}},
volume = {6},
journal = {The Annals of Statistics},
number = {4},
publisher = {Institute of Mathematical Statistics},
pages = {882 -- 894},
year = {1978},
doi = {10.1214/aos/1176344261}
}

@article{bickel1981someAsymptotic,
author = {Peter J. Bickel and David A. Freedman},
title = {{Some asymptotic theory for the bootstrap}},
volume = {9},
journal = {The Annals of Statistics},
number = {6},
publisher = {Institute of Mathematical Statistics},
pages = {1196 -- 1217},
year = {1981},
doi = {10.1214/aos/1176345637}
}

@article{Weed2019sharpAsymptotic,
author = {Weed, Jonathan and Bach, Francis},
title = {{Sharp asymptotic and finite-sample rates of convergence of empirical measures in Wasserstein distance}},
year = {2019},
journal = {Bernoulli},
volume = {25},
number = {4 A},
pages = {2620 – 2648},
doi = {10.3150/18-BEJ1065}
}

@article{Bonneel2015slicedAndRadon,
author={Bonneel, Nicolas
and Rabin, Julien
and Peyr{\'e}, Gabriel
and Pfister, Hanspeter},
title={{Sliced and Radon Wasserstein barycenters of measures}},
journal={Journal of Mathematical Imaging and Vision},
year={2015},
day={01},
volume={51},
number={1},
pages={22-45},
issn={1573-7683},
doi={10.1007/s10851-014-0506-3}
}

@article{delBarrio2023improvedEntropic,
author = {del Barrio, Eustasio and Sanz, Alberto Gonz\'{a}lez and Loubes, Jean-Michel and Niles-Weed, Jonathan},
title = {{An improved central limit theorem and fast convergence rates for entropic transportation costs}},
journal = {SIAM Journal on Mathematics of Data Science},
volume = {5},
number = {3},
pages = {639-669},
year = {2023},
doi = {10.1137/22M149260X}
}

@article{hundrieser2024unifying,
author = {Shayan Hundrieser and Marcel Klatt and Axel Munk and Thomas Staudt},
title = {{A unifying approach to distributional limits for empirical optimal transport}},
volume = {30},
journal = {Bernoulli},
number = {4},
publisher = {Bernoulli Society for Mathematical Statistics and Probability},
pages = {2846 -- 2877},
year = {2024},
doi = {10.3150/23-BEJ1697}
}

@article{Sen1977Ustatistics,
author = {Pranab Kumar Sen},
title = {{Almost sure convergence of generalized $U$-statistics}},
volume = {5},
journal = {The Annals of Probability},
number = {2},
publisher = {Institute of Mathematical Statistics},
pages = {287 -- 290},
year = {1977},
doi = {10.1214/aop/1176995853}}

@article{csorgo1990distributions,
     author = {Cs\"org\"o, Mikl\'os and Horv\'ath, Lajos},
     title = {{On the distributions of $L^p$ norms of weighted quantile processes}},
     journal = {Annales de l'I.H.P. Probabilit\'es et statistiques},
     pages = {65--85},
     publisher = {Gauthier-Villars},
     volume = {26},
     number = {1},
     year = {1990},
     mrnumber = {1075439},
     zbl = {0699.62045},
     language = {en}
}

@article{weed2022estimation,
author = {Jonathan Niles-Weed and Philippe Rigollet},
title = {{Estimation of Wasserstein distances in the Spiked Transport Model}},
volume = {28},
journal = {Bernoulli},
number = {4},
publisher = {Bernoulli Society for Mathematical Statistics and Probability},
pages = {2663 -- 2688},
year = {2022},
doi = {10.3150/21-BEJ1433}
}

@article{xu2022central,
  title={{Central limit theorem for the sliced 1-Wasserstein distance and the max-sliced 1-Wasserstein distance}},
  author={Xu, Xianliang and Huang, Zhongyi},
  journal={arXiv preprint},
  year={2022},
  doi = {
https://doi.org/10.48550/arXiv.2205.14624}
}

@article{gangbo1996geometry,

author = {Wilfrid Gangbo and Robert J. McCann},
title = {{The geometry of optimal transportation}},
volume = {177},
journal = {Acta Mathematica},
number = {2},
publisher = {Institut Mittag-Leffler},
pages = {113 -- 161},
year = {1996},
doi = {10.1007/BF02392620}
}

@inproceedings{mena2019statisticalBounds,
 author = {Mena, Gonzalo and Niles-Weed, Jonathan},
 booktitle = {Advances in Neural Information Processing Systems},
 publisher = {Curran Associates, Inc.},
 title = {{Statistical bounds for entropic optimal transport: sample complexity and the central limit theorem}},
 volume = {32},
 year = {2019}
}

@inproceedings{Cuturi2013SinkhornDL,
 author = {Cuturi, Marco},
 booktitle = {Advances in Neural Information Processing Systems},
 pages = {},
 publisher = {Curran Associates, Inc.},
 title = {{Sinkhorn Distances: Lightspeed computation of optimal transport}},
 volume = {26},
 year = {2013}
}

@book{villani2003Topics,
author = {Villani, Cédric},
address = {Providence, R.I},
booktitle = {Topics in optimal transportation},
isbn = {082183312X},
language = {eng},
lccn = {2003040350},
publisher = {American Mathematical Society},
series = {Graduate studies in mathematics},
title = {{Topics in optimal transportation}},
year = {2003},
}

@book {villani2008optimal,
    AUTHOR = {Villani, C\'edric},
     TITLE = {Optimal transport. Old and new.},
 PUBLISHER = {Springer},
      YEAR = {2009}
}

@book{Korolyuk2013theoryUstatistic,
  title={{Theory of U-statistics}},
  author={Korolyuk, Vladimir S. and Borovskich, Yu V.},
  volume={273},
  year={2013},
  publisher={Springer Science \& Business Media},
doi = {https://doi.org/10.1007/978-94-017-3515-5}
}

@book{shorack2009empirical,
  title={{Empirical processes with applications to statistics}},
  author={Shorack, Galen R. and Wellner, Jon A.},
  isbn={9780898716849},
  lccn={2009025143},
  series={Classics in Applied Mathematics},
  year={2009},
  publisher={Society for Industrial and Applied Mathematics}
}

@article{sommerfeld2017inferenceFiniteSpaces,
    author = {Sommerfeld, Max and Munk, Axel},
    title = {{Inference for empirical Wasserstein distances on finite spaces}},
    journal = {Journal of the Royal Statistical Society Series B: Statistical Methodology},
    volume = {80},
    number = {1},
    pages = {219-238},
    year = {2017},
    issn = {1369-7412},
    doi = {10.1111/rssb.12236}
}

@article{tameling2019empirical,
author = {Carla Tameling and Max Sommerfeld and Axel Munk},
title = {{Empirical optimal transport on countable metric spaces: Distributional limits and statistical applications}},
volume = {29},
journal = {The Annals of Applied Probability},
number = {5},
publisher = {Institute of Mathematical Statistics},
pages = {2744 -- 2781},
year = {2019},
doi = {10.1214/19-AAP1463}
}

@article{delBarrio1999,
author = {Eustasio del Barrio and Evarist Gin{\'e} and Carlos Matr{\'a}n},
title = {{Central limit theorems for the Wasserstein distance between the empirical and the true distributions}},
volume = {27},
journal = {The Annals of Probability},
number = {2},
publisher = {Institute of Mathematical Statistics},
pages = {1009 -- 1071},
year = {1999},
doi = {10.1214/aop/1022677394}
}

@article{delBarrio2005asyptotics,
author = {Eustasio del Barrio and Evarist Gin{\'e} and Frederic Utzet},
title = {{Asymptotics for $L^2$ functionals of the empirical quantile process, with applications to tests of fit based on weighted Wasserstein distances}},
volume = {11},
journal = {Bernoulli},
number = {1},
publisher = {Bernoulli Society for Mathematical Statistics and Probability},
pages = {131 -- 189},
year = {2005},
doi = {10.3150/bj/1110228245}
}

@article{BonneelDigne2023,
author = {Bonneel, Nicolas and Digne, Julie},
title = {A survey of Optimal Transport for Computer Graphics and Computer Vision},
journal = {Computer Graphics Forum},
volume = {42},
pages = {439-460},
year = {2023}
}

@inproceedings{gordaliza2019obtaining,
  title={Obtaining fairness using optimal transport theory},
  author={Gordaliza, Paula and Del Barrio, Eustasio and Fabrice, Gamboa and Loubes, Jean-Michel},
  booktitle={International conference on machine learning},
  pages={2357--2365},
  year={2019},
  organization={PMLR}
}

@article{chang2022unified,
  title={Unified optimal transport framework for universal domain adaptation},
  author={Chang, Wanxing and Shi, Ye and Tuan, Hoang and Wang, Jingya},
  journal={Advances in Neural Information Processing Systems},
  volume={35},
  pages={29512-29524},
  year={2022}
}

@article{deLara.et.al.JMLR.2024,
  title={Transport-based counterfactual models},
  author={de Lara, Lucas and Gonz{\'a}lez-Sanz, Alberto and Asher, Nicholas and Risser, Laurent and Loubes, Jean-Michel},
  journal={Journal of Machine Learning Research},
  volume={25},
  pages={1--59},
  year={2024}
}

@article{freulon2023cytopt,
  title={CytOpT: Optimal transport with domain adaptation for interpreting flow cytometry data},
  author={Freulon, Paul and Bigot, J{\'e}r{\'e}mie and Hejblum, Boris P},
  journal={Ann. Appl. Statist.},
  volume={17},
  pages={1086--1104},
  year={2023}
}

@article{del2020optimalflow,
  title={optimalFlow: optimal transport approach to flow cytometry gating and population matching},
  author={del Barrio, Eustasio and Inouzhe, Hristo and Loubes, Jean-Michel and Matr{\'a}n, Carlos and Mayo-{\'I}scar, Agust{\'\i}n},
  journal={BMC bioinformatics},
  volume={21},
  pages={1--25},
  year={2020},
}

@article{pooladian2024pluginestimationschrodingerbridges,
  title={{Plug-in Estimation of Schr{\"o}dinger Bridges}},
  author={Pooladian, Aram-Alexandre and Niles-Weed, Jonathan},
  journal={SIAM Journal on Mathematics of Data Science},
  volume={7},
  number={3},
  pages={1315--1336},
  year={2025},
  publisher={SIAM}
}

@inproceedings{tolstikhin2018wasserstein,
  title={Wasserstein Auto-Encoders},
  author={Tolstikhin, Ilya and Bousquet, Olivier and Gelly, Sylvain and Schoelkopf, Bernhard},
  booktitle={International Conference on Learning Representations},
  year={2018}
}

@inproceedings{arjovsky2017wasserstein,
  title={Wasserstein generative adversarial networks},
  author={Arjovsky, Martin and Chintala, Soumith and Bottou, L{\'e}on},
  booktitle={International conference on machine learning},
  pages={214--223},
  year={2017},
  organization={PMLR}
}

@article{COTFNT, 
year = {2019}, 
volume = {11}, 
journal = {Foundations and Trends in Machine Learning}, 
title = {Computational Optimal Transport}, number = {5-6}, 
pages = {355--607},
author = {Gabriel Peyr\'e and Marco Cuturi}
}

@article{alvarez2011uniqueness,
     author = {\'Alvarez-Esteban, P. C. and del Barrio, E. and Cuesta-Albertos, J. A. and Matr\'an, C.},
     title = {Uniqueness and approximate computation of optimal incomplete transportation plans},
     journal = {Annales de l'I.H.P. Probabilit\'es et statistiques},
     pages = {358--375},
     publisher = {Gauthier-Villars},
     volume = {47},
     number = {2},
     year = {2011},
     doi = {10.1214/09-AIHP354},
     mrnumber = {2814414},
     zbl = {1215.49042},
     language = {en},
     url = {https://www.numdam.org/articles/10.1214/09-AIHP354/}
}

\appendix
\section*{APPENDIX}

All results presented in this work apply to both the one-sample and two-sample settings, with the proofs being largely similar in most cases. The following sections of the appendix provide detailed proofs for the one-sample setting and highlight the key modifications required to extend these arguments to the two-sample setting.

\section{Proofs of Section \ref{sect:OT_potentials}}\label{appendix:OT_potentials}

\subsection{Preliminary results}

\begin{lemma} \label{lemma:support}
    Let $P\in\mathcal P(\mathbb R^d)$ be an absolutely continuous probability with negligible boundary and $\operatorname{int}(\supp P)$ connected. Then, for every $\theta\in\mathbb S^{d-1}$, $P^\theta$ is also absolutely continuous, has negligible boundary, and $\operatorname{int}(\supp{P^\theta})$ is connected.
\end{lemma}
\begin{proof}
First, we prove the inclusions $\operatorname{Pr}_\theta\bigl( \operatorname{int}(\supp P)\bigr) \underset{(i)}{\subset}  \operatorname{int}(\supp{P^\theta})  \underset{(ii)}{\subset} \operatorname{cl}\bigl( \operatorname{Pr}_\theta\bigl( \operatorname{int}(\supp P)\bigr)\bigr) \ .$
\begin{enumerate}[label=(\roman*)]
    \item Let $s\in\operatorname{Pr}_\theta\bigl( \operatorname{int}(\supp P)\bigr)$. There exist $x\in \operatorname{int}(\supp P)$ and $\varepsilon>0$ such that $s=\langle \theta,x\rangle$ and $B(x,\varepsilon)\subset \supp P$. Then, for every $s'=s+t\in B(s,\varepsilon)$, $x'= x+\theta t \in B(x,\varepsilon)$.  Therefore, $x'\in \supp P$. Thus, for every $r>0$, 
    \begin{equation*}
        P^\theta\bigl(B(s',r)\bigr)= P\bigl( \{x\in\mathbb R^d: \langle \theta,x\rangle \in B(s',r)\}\bigr)\geq P\bigl( B(x',r)\bigr) >0 \ .
    \end{equation*}
    Therefore, $s'\in\supp{P^\theta}$, and it follows that $B(s,\varepsilon)\subset \supp{P^\theta}$, which implies that $s\in  \operatorname{int}(\supp{P^\theta})$.
    
    \item Let $s\in  \operatorname{int}(\supp{P^\theta})$. For each $n\in\mathbb N$, since $P$ is absolutely continuous with  negligible boundary, $P(\operatorname{int}(\supp P))=1$, which implies 
    \begin{align*}
         0<P^\theta\bigl(B(s,1/n)\bigr) &=P\bigl( \{x\in\mathbb R^d: \langle \theta,x\rangle \in B(s,1/n)\}\bigr)=P\bigl( \{x\in \operatorname{int}(\supp P): \langle \theta,x\rangle \in B(s,1/n)\}\bigr) \ .
        \end{align*}
        Thus, we can define a sequence of values $x_n \in \operatorname{int}(\supp P)$ such that $ \langle \theta,x_n\rangle \rightarrow s$. Therefore, $s$ belongs to the closure of $\operatorname{Pr}_\theta\bigl( \operatorname{int}(\supp P)\bigr)$.
\end{enumerate}
Now, since $\operatorname{int}(\supp P)$ is connected and $\operatorname{Pr}_\theta$ is open and continuous, it follows that $\operatorname{Pr}_\theta\bigl( \operatorname{int}(\supp P)\bigr)$ is an open interval $I^\theta$, which, together with the previous inclusions, allows us to conclude
\begin{equation}    \label{eq:supports}
    \operatorname{int}(\supp{P^\theta}) = \operatorname{Pr}_\theta\bigl( \operatorname{int}(\supp P)\bigr)= I^\theta \ .
\end{equation}  
Therefore, $P^\theta$ has negligible boundary and $\operatorname{int}(\supp{P^\theta})$ is connected. The absolute continuity is straightforward.
\end{proof}

The following lemma summarizes key regularity properties of $c$-concave optimal transport potentials used in the proof of Proposition \ref{prop:regularidad_potenciales}. Assertions \textit{(i)}, \textit{(ii)} and \textit{(iii)} are just an adaptation of Theorem 3.3 and Proposition C.4 in \cite{gangbo1996geometry} for $d=1$, but these properties hold more generally.

\begin{lemma}\label{lemma:c-convave}
    Let $p>1$, $P,Q\in\mathcal P_p(\mathbb R)$ and $c=c_p$. Denote by $\phi$ any $c$-concave optimal transport potential from $P$ to $Q$. Then, 
    \begin{enumerate}[label=(\roman*)]
        \item There exists a closed interval $I$ such that $\operatorname{int}(I)\subset \operatorname{dom} \phi = \{x\in\mathbb R : \phi(x)>-\infty\} \subset I$.
        \item $\phi$ is locally Lipschitz on $\operatorname{int}(I)$, and, therefore, continuous.
        \item For every $s\in \operatorname{int}(I)$, $\partial^c\phi(s)\neq \emptyset$.
        \item Given $s\in\mathbb R$, if $\partial^c\phi(s)\neq \emptyset$, then $\partial^c\phi(s)$ is a single point or a closed interval.
    \end{enumerate}
\end{lemma}

\begin{proof}
    The only assertion that requires to be proved is \textit{(iv)}. Assume that $\partial^c\phi(s)$ does not consist of a single point. To prove that $\partial^c\phi(s)$ is an interval, we need to see that if $t_1,t_2\in \partial^c\phi(s)$, and $t\in [t_1,t_2]$, then $t\in\partial^c\phi(s)$. This 
 follows from the definition \eqref{eq:defn_superdifferential} and the inequality 
\begin{equation*}
    \min \Big\{  |z-t_1|^p-|x-t_1|^p ,  |z-t_2|^p-|x-t_2|^p \Bigr\} \leq  |z-t|^p-|x-t|^p \ ,
\end{equation*}
which follows from the fact that $t\mapsto  |z-t|^p-|x-t|^p$ is monotone increasing if $z\leq x$ and decreasing if $z\geq x$. Thus, $\partial^c\phi(s)$ is an interval. To see that it is closed, let $t\in \operatorname{cl}(\partial^c\phi(s))$, and consider a sequence $\{t_n\}_{n=1}^\infty \subset \partial^c\phi(s)$, such that $t_n\rightarrow t$. Then, for every $z\in \mathbb R$, $\phi(z) \leq \phi(s) + \bigl( |z-t_n|^p-|x-t_n|^p\bigr)$.
By continuity, we conclude that $t\in \partial^c \phi(x)$.
\end{proof}

With the notation of the previous  lemma, it follows that $\supp{P}\subset I$. Otherwise, if $s\in \supp P \setminus I$, by definition of $\supp P$, for every $\varepsilon>0$, $P(B(s,\varepsilon))>0$. Since $I$ is closed, $\mathbb R\setminus I$ is open, and there exists $\varepsilon>0$ such that $B(s,\varepsilon)\subset \mathbb R\setminus I \subset \mathbb R\setminus \operatorname{dom}(\phi)$. This implies that $\int \phi(s)dP(s)$ is not finite, which contradicts the optimality of $\phi$.

Therefore, $\operatorname{int}(\supp{P})\subset \operatorname{int}(I)$. If $P$ is absolutely continuous, $\operatorname{int}(\supp{P})$ defines a set of probability one, and Lemma \ref{lemma:c-convave} allows us to guarantee favorable properties on this set. Now, if $P$ assigns positive probability to $s\in\mathbb R$, by the optimality, $\phi(s)\in \operatorname{dom}(f)$. In this case, even if $s\in I\setminus \operatorname{int}(I)$, we can ensure $\partial^c\phi(s)\neq \emptyset$. To check this, let $\pi\in\Pi(P,Q)$ be any optimal transport plan such that $\supp\pi\subset \partial^c\phi$. Since $P(s)>0$, we know that $\pi(\{s\}\times \mathbb R)>0$. Therefore, there exists $A\subset R$ such that $\{s\}\times A\subset \supp\pi\subset \partial^c\phi$, which implies $A\subset \partial^c \phi(s)$.

\subsection{Proof of Proposition \ref{prop:regularidad_potenciales}}
\noindent
We restrict our attention to the subset of probability one $\Omega_0\subset \Omega$  verifying that for every $\omega\in\Omega_0$,
$X_1^\omega\in \operatorname{int}(\supp{P})$ 
and 
$P_n^\omega \rightsquigarrow P$. 
Thus, it suffices to show all the assertions for a sequence $\{X_n^\omega\}_{n=1}^\infty$, with $\omega\in\Omega_0$.  For simplicity, we omit the superscript $\omega$.
Note that some parts of the proof closely follow the arguments of Lemma B.2 in \cite{Goldfel2024statisticalInference}.
 \\

\noindent
{(i)} For each $\theta\in\mathbb S^{d-1}$, by Lemma \ref{lemma:support}, $P^\theta$ is absolutely continuous with negligible boundary and $\operatorname{int}(\supp{P^\theta})$ connected. By the comments following Lemma \ref{lemma:c-convave}, $\phi^\theta$ is finite on $\operatorname{int}(\supp{P^\theta})$. Since $X_1\in \operatorname{int}(\supp{P})$, by \eqref{eq:supports}, $\langle\theta,X_1\rangle \in \operatorname{int}(\supp{P^\theta})$. Therefore, $\phi^\theta(\langle \theta,X_1\rangle)\in\mathbb  R$ and $\phi^\theta_{X_1}(s) = \phi^\theta(s)-\phi^\theta(\langle\theta,X_1\rangle)$ is well defined, it is a $c$-concave optimal transport potential and verifies $\phi^\theta_{X_1}(\langle \theta,X_1\rangle)=0
$. Moreover, if $\Tilde \phi_\theta$ is any other $c$-concave optimal transport potential verifying $\Tilde \phi^\theta(\langle \theta,X_1\rangle)=0$, by Corollary 2.7 in \cite{delBarrio2024CLTgeneral}, there exists $a^\theta \in \mathbb R$ such that  $\phi^\theta_{X_1} = \Tilde \phi^\theta + a^\theta$ on  $\operatorname{int}(\supp{P^\theta})$. Since both potentials vanish at $\langle \theta,X_1\rangle$, it follows that $a^\theta=0$. \\

\noindent
{(ii)}  Let $x\in \operatorname{int}(\supp{P})$, and $\theta_n \rightarrow \theta$. It is easy to see that $P^{\theta_n}\rightsquigarrow P^{\theta}$. Thus, the first part of Theorem 3.4 in \cite{delBarrio2024CLTgeneral} implies that there exist constants $a^{\theta_n}\in\mathbb R$ (which depend on $\omega\in\Omega_0$) such that $\phi_{X_1}^{\theta_n} - a^{\theta_n} \rightarrow \phi_{X_1}^\theta$ in the sense of uniform convergence on the compact sets of $\operatorname{int}(\supp{P^\theta})$.  By the triangle inequality,
\begin{align*}
    |a^{\theta_n}| \leq  |(\phi_{X_1}^{\theta_n}(\langle\theta_n,X_1\rangle)-a^{\theta_n})-\phi_{X_1}^\theta(\langle\theta_n,X_1\rangle)| + |\phi_{X_1}^{\theta_n}(\langle\theta_n,X_1\rangle)-\phi_{X_1}^\theta(\langle\theta_n,X_1\rangle)|  \ .
\end{align*}
The first term converges to zero by the uniform convergence on compact sets of $\operatorname{int}(\supp{P^\theta})$ and the fact that $\langle \theta_n, X_1 \rangle \to \langle \theta, X_1 \rangle \in \operatorname{int}(\supp{P^\theta})$, whereas the second term converges to zero by  $\phi_{X_1}^{\theta_n}(\langle\theta_n,X_1\rangle)=\phi_{X_1}^\theta(\langle\theta,X_1\rangle)=0$ and the continuity of $\phi_{X_1}^\theta$ on $\operatorname{int}(\supp{P^\theta})$, guaranteed by  Lemma \ref{lemma:c-convave}. Therefore, $a^{\theta_n}\rightarrow 0$, and we can conclude that 
\begin{align*}
|\phi_{X_1}^{\theta_n}-\phi_{X_1}^\theta| \leq  |(\phi_{X_1}^{\theta_n}-a^{\theta_n})-\phi_{X_1}^\theta|+|a^{\theta_n}|\rightarrow 0
\end{align*}
uniformly on the compact sets of $\operatorname{int}(\supp{P^\theta})$. Finally, for any $x \in \operatorname{int}(\supp{P})$, since $\langle \theta, x \rangle \in \operatorname{int}(\supp{P^\theta})$ by \eqref{eq:supports}, we obtain
\begin{align}\label{eq:prop2.1(2)}
  &|\phi^{\theta_n}_{X_1}(\langle\theta_n,x\rangle) - \phi_{X_1}^\theta(\langle\theta,x\rangle) | \leq |\phi^{\theta_n}_{X_1}(\langle\theta_n,x\rangle) - \phi_{X_1}^\theta(\langle\theta_n,x\rangle) |+ |\phi^{\theta}_{X_1}(\langle\theta_n,x\rangle) - \phi_{X_1}^\theta(\langle\theta,x\rangle) | \rightarrow 0 
\end{align} 
where the first term vanishes again by uniform convergence on the compact sets of $\operatorname{int}(\supp{P^\theta})$, and the second by continuity of $\phi^\theta_{X_1}$ on $\operatorname{int}(\supp{P^\theta})$.

\noindent
{(iii)} By the discussion after Lemma \ref{lemma:c-convave},  $\phi_n^\theta(\langle\theta,X_1\rangle)<\infty$. Therefore, $\phi_{n,X_1}^\theta(s)=\phi_{n}^\theta(s)- \phi_n^\theta(\langle\theta,X_1\rangle)$ is well defined, it is a $c$-concave optimal transport potential and verifies $\phi^\theta_{n,X_1}(\langle \theta,X_1\rangle)=0
$. \\
 

\noindent
{(iv)} 
For each $\theta\in\mathbb S^{d-1}$, the continuous mapping theorem implies $P_n^\theta \rightsquigarrow P^\theta$. As in the proof of (2), there exists constants $a_n^\theta\in\mathbb R$ such that $\phi_{X_1}^{\theta_n} - a_n^\theta \rightarrow \phi_{X_1}^\theta$ in the sense of uniform convergence on the compact sets of $\operatorname{int}(\supp{P^\theta})$. Since  $\phi^{\theta}_{n,X_1}(\langle \theta,X_1\rangle)=\phi^{\theta}_{X_1}(\langle \theta,X_1\rangle)=0$, it follows that  $a_n^\theta\rightarrow 0$ and, therefore, $\phi^{\theta}_{n,X_1}(t) \rightarrow\phi^\theta_{X_1}(t)$
uniformly on the compact sets of $\operatorname{int}(\supp{P^\theta})$. The second assertion follows directly from the second part of Theorem 3.4 in \cite{delBarrio2024CLTgeneral}.\\


\noindent
{(v)} 
Assume that the property is not verified. Then, there exists a subsequence $\{n_k\}_{k=1}^\infty$ and values $\{\theta_{n_k}\}_{k=1}^\infty$ such that $    |\phi_{n_k,X_1}^{\theta_{n_k}}( \langle \theta_{n_k}, x\rangle)| \rightarrow \infty$. 
Since $\mathbb S^{d-1}$ is compact, we can extract a further subsequence $\theta_{n_{k_m}}\rightarrow \theta \in \mathbb S^{d-1}$. For this subsequence,  $P^{\theta_{n_{k_m}}}_{n_{k_m}} \rightsquigarrow P^\theta$ and $Q^{\theta_{n_{k_m}}} \rightsquigarrow Q^\theta$. By the same argument used in the proof of (2), it follows that $\phi_{n_{k_m},X_1}^{\theta_{n_{k_m}}}(t)$ converges to $\phi^\theta_{X_1}(t)$ uniformly  on the compact sets of $\operatorname{int}(\operatorname{supp}(P^\theta))$, and reasoning as in \eqref{eq:prop2.1(2)}, we can conclude 
$
    \phi_{n_{k_m},X_1}^{\theta_{n_{k_m}}}( \langle \theta_{n_{k_m}}, x\rangle) \longrightarrow \phi^{\theta}_{X_1}( \langle \theta, x\rangle) < \infty \ ,
$
which contradicts the hypothesis.\\

\noindent
{(vi)} For each $m\in\mathbb N$, the existence of $\phi_{n,m,X_1}$ follows from Property (3) with $Q=Q_m$. Theorem 3.4 in \cite{delBarrio2024CLTgeneral} also applies to the two-sample setting. Therefore, the same proof of (4), (5) and (6) extends to this case in the set of $\mathbb P$-probability one $\Omega_0 \cap \{\omega: Q_m\rightsquigarrow Q\}$.
\qed

\subsection{Proof of Proposition \ref{prop:integrability_potentials}}

\noindent

First, we prove \eqref{eq:int_potentials_one_sample}.  By Proposition \ref{prop:regularidad_potenciales}, $\phi^\theta_{X_1}$ and $\phi^\theta_{n,X_1}$ are well-defined $c$-concave optimal transport potentials with reference point $X_1$, up to a set of $\mathbb P$-probability zero. To prove integrability properties, we can assume without loss of generality that it is verified for every point. Then, by the $c$-concavity,
\begin{equation}\label{eq:prop_Villani}
    |s-t|^p=\phi_{n,X_1}^\theta(s)+(\phi_{n,X_1}^\theta)^c(t) \Longleftrightarrow t \in \partial^c \phi_{n,X_1}^\theta(s) \Longleftrightarrow s \in \partial^c (\phi_{n,X_1}^\theta)^c(t) \ .
\end{equation}
Let $t_{n,X_1}^\theta$ be any element of $\partial^c \phi_{n,X_1}^\theta\left(\langle\theta,X_1\rangle\right)$, which is non-empty by the comments after Lemma \ref{lemma:c-convave}. The first part of the argument is valid for any choice of $t_{n,X_1}^\theta\in \partial^c \phi_{n,X_1}^\theta\left(\langle\theta,X_1\rangle\right)$, although to conclude, we will consider a particular choice in \eqref{eq;adequate_subgradient}.  
By definition of the $c$-superdifferential, for every $s\in\mathbb R$,
\begin{equation}\label{eq:subdiff}
    \phi_{n,X_1}^\theta(s) \leq \phi_{n,X_1}^\theta\left(\langle \theta,X_1 \rangle\right)+\Bigl( |s-t_{n,X_1}^\theta|^p - |\langle \theta,X_1 \rangle- t_{n,X_1}^\theta|^p\Bigr) \ .
\end{equation}
Since $t_{n,X_1}^\theta\in \partial^c \phi_{n,X_1}^\theta\left(\langle \theta,X_1 \rangle\right)$, it follows by   \eqref{eq:prop_Villani} that $\langle \theta,X_1 \rangle \in \partial^c (\phi_{n,X_1}^\theta)^c(t_{n,X_1}^\theta)$. Therefore, for every $t\in\mathbb R$,
\begin{equation}\label{eq:subdiff_conjugada}
   (\phi_{n,X_1}^\theta)^c(t) \leq (\phi_{n,X_1}^\theta)^c(t_{n,X_1}^\theta)+\Bigl( |t-\langle \theta,X_1 \rangle|^p - | t_{n,X_1}^\theta-\langle \theta,X_1 \rangle|^p\Bigr) \ .
\end{equation}
By \eqref{eq:prop_Villani} and \eqref{eq:subdiff_conjugada}, it follows that for every $(s,t) \in \partial^c \phi_{n,X_1}^\theta$,
\begin{align} \label{eq:cota_phi_n}
    \phi_{n,X_1}^\theta(s) & =|s-t|^p  - (\phi_{n,X_1}^\theta)^c(t) \geq\phi_{n,X_1}^\theta\left(\langle \theta,X_1 \rangle\right)+ |s-t|^p - |t-\langle \theta,X_1 \rangle|^p
\end{align}
For every $(s,t) \in \partial^c \phi_{n,X_1}^\theta$, combining \eqref{eq:subdiff} and \eqref{eq:cota_phi_n} and using inequality \eqref{eq:elemental_ineq_p},
\begin{align}\label{eq:cota_buena}
&\big|\phi_{n,X_1}^\theta(s) - \phi_{n,X_1}^\theta(\langle \theta,X_1 \rangle) \big|  \leq 2C_p \Bigl( |s|^p + |t|^p + |\langle \theta,X_1 \rangle|^p +|t_{n,X_1}^\theta|^p \Bigr) \ .
\end{align}
Similarly, for every $(s,t) \in \partial^c \phi_{n,X_1}^\theta$,
\begin{align}\label{eq:cota_buena_conjugate}
&\big|(\phi_{n,X_1}^\theta)^c(t) - (\phi_{n,X_1}^\theta)^c(t_{n,X_1}^\theta) \big|  \leq 2C_p \Bigl( |s|^p + |t|^p + |\langle \theta,X_1 \rangle|^p +|t_{n,X_1}^\theta|^p \Bigr) \ .
\end{align}
But optimality implies that there exists an optimal coupling $\pi_n^\theta$ between $P_n^\theta$ and $Q^\theta$ such that $\supp{\pi_n^\theta} \subset \partial^c \phi_{n,X_1}^\theta$. Therefore, \eqref{eq:cota_buena} is verified for every $(s,t)\in \supp{\pi_n^\theta}$. Given that $\phi_{n,X_1}^\theta(\langle\theta,X_1\rangle)=0$ and $\phi_{n,X_1}^\theta(\langle\theta,X_1\rangle)+(\phi_{n,X_1}^\theta)^c(t_{n,X_1}^\theta)=|\langle\theta,X_1\rangle -t_{n,X_1}^\theta|^p$ it follows that both integrals 
\begin{align*}
    \int_{\mathbb R^d} \big|\phi_{n,X_1}^\theta(\langle \theta,x\rangle) \big|^{q} dP_n(x)  &=  \int_{\mathbb R\times \mathbb R} \big|\phi_{n,X_1}^\theta(s) \big|^{q} d\pi_n^\theta(s,t) \\
    \int_{\mathbb R^d} \big|(\phi_{n,X_1}^\theta)^c(\langle \theta,y\rangle) \big|^{q} dQ(y)  &=  \int_{\mathbb R\times \mathbb R} \big|(\phi_{n,X_1}^\theta)^c(t) \big|^{q} d\pi_n^\theta(s,t)
\end{align*}
are bounded, up to a constant, by
\begin{align*}
 \int_{\mathbb R} |s|^{pq} dP_n^\theta(s) + \int_{\mathbb R} |t|^{pq} dQ^\theta(t) +|\langle \theta,X_1 \rangle|^{pq} +|t_{n,X_1}^\theta|^{pq} \ .
\end{align*}
Therefore, we only need to uniformly bound the expected value of the integral over $\mathbb S^{d-1}$ of each term. 
\begin{itemize}
    \item $\mathbb E \Bigl( \int_{\mathbb S^{d-1}} \int_{\mathbb R} |s|^{pq} dP_n^\theta(s) d\sigma(\theta) \Bigr)   \leq\int_{\mathbb R^d} \|x\|^{pq} dP(x) \ .$\\
\item $\mathbb E \Bigl( \int_{\mathbb S^{d-1}} \int_{\mathbb R} |t|^{pq} dQ^\theta(t) d\sigma(\theta) \Bigr) \leq \int_{\mathbb R^d} \|y\|^{pq} dQ(y) \ .$\\
\item  $\mathbb E \Bigl( \int_{\mathbb S^{d-1}} |\langle \theta,X_1 \rangle|^{pq} d\sigma(\theta) \Bigr) \leq \int_{\mathbb R^d} \|x\|^{pq} dP(x) \ .$\\
\item  For the last term, we need to consider a particular choice of the elements $t_{n,X_1}^\theta$. Note that the $c$-superdifferential does not depend on the particular choice of the constant. Thus, $\partial^c \phi_{n,X_1}^\theta(\langle \theta,x\rangle) = \partial^c \phi_{n}^\theta(\langle \theta,x\rangle)$ and we define
for each $i=1,\ldots,n$ and $\theta\in\mathbb S^{d-1}$ the random variable 
\begin{equation}\label{eq;adequate_subgradient}
         t_{n,X_i}^\theta = \underset{t\in \partial^c \phi_{n}^\theta(\langle\theta,X_i\rangle)}{\operatorname{argmin}} |t| \ ,    \end{equation}
     which is well defined since the $c$-superdifferential of a $c$-concave function is a closed (possibly infinite) interval, by Lemma \ref{lemma:c-convave}.  Note that, by definition,
$
          t_{n,X_1}^\theta\overset{d}{=}\ldots \overset{d}{=} t_{n,X_n}^\theta \ .
  $
     Then, given an optimal transport plan $\pi_n^\theta$ verifying  $\operatorname{supp}(\pi_n^\theta)\subset \partial^c \phi_n^\theta$, 
     \begin{align}
        \int_{\mathbb R^d} \|y\|^{pq} dQ(y) &\geq \int_{\mathbb R^d}  \left| \langle \theta,y\rangle \right|^{pq} dQ(y)=\int_{\mathbb R\times \mathbb R}  \left| t \right|^{pq} d\pi_n^\theta(s,t) = \sum_{i=1}^n \int_{\{ s=\langle \theta,X_i\rangle, t\in \partial^c\phi_n^\theta(s)\}}  \left| t \right|^{pq} d\pi_n^\theta(s,t)\notag\\
        &\geq \sum_{i=1}^n \int_{\{ s=\langle \theta,X_i\rangle, t\in \partial^c\phi_n^\theta(s)\}}  \bigl| t_{n,X_i}^\theta \bigr|^{pq} d\pi_n^\theta(s,t) = \frac{1}{n}\sum_{i=1}^n \bigl| t_{n,X_i}^\theta \bigr|^{pq}  \ , \label{eq:ineq_measurability_issues}
     \end{align}
which allow us to conclude
    $
         \mathbb E(| t_{n,X_1}^\theta|^{pq} ) =  \mathbb E(\tfrac{1}{n}\sum_{i=1}^n | t_{n,X_i}^\theta |^{pq}  ) \leq  \int_{\mathbb R^d} \|y\|^{pq} dQ(y) \ .
    $
\end{itemize}

\noindent
The proof of \eqref{eq:int_potentials_limit} follows similarly. For each $\theta\in\mathbb S^{d-1}$, by Theorem 1.2 in \cite{gangbo1996geometry},  $\partial^c \phi_\theta(s) = \{ \nabla^c \phi_\theta(s) \}$ $P^\theta$-a.s. and $\nabla^c \phi^\theta$ is an optimal transport map from $P^\theta$ to $Q^\theta$. Therefore,  $t_{X_1}^\theta := \nabla^c \phi^\theta(\langle \theta, X_1 \rangle)$
is well defined $(\mathbb P\times\sigma)$-a.s. and 
\begin{align*}
    \mathbb E \Bigl( \int_{\mathbb S^{d-1}} |t_{X_1}^\theta|^{pq} d\sigma(\theta)\Bigr)    &  = \int_{\mathbb S^{d-1}}  \Bigl( \int_{\mathbb R} |\nabla^c \phi_\theta(s)|^{pq} dP^\theta(s)\Bigr)d\sigma(\theta) 
     \leq \int_{\mathbb R^d} \|y\|^{pq} dQ(y)
\end{align*}  

\noindent
\textit{Two-sample setting:} \eqref{eq:int_potentials_two_sample} follows from the same reasoning, using that $\mathbb E \bigl( \int_{\mathbb R} \|y\|^{pq} dQ_m(y)\bigr)  = \int_{\mathbb R} \|y\|^{pq} dQ(y) $.
\qed

\section{Proofs of Section \ref{sect:CLT_general}}\label{appendixA}

\subsection{Proof of Proposition \ref{teor:primal}}

Define the random variable $Z_n= \sw_p^p(P_n,Q)$, where $P_n= \frac{1}{n}\sum_{i=1}^n \delta_{X_i}$, with $X_1,\ldots,X_n$ i.i.d. random variables with distribution $P$. Consider now  $Z_n'= \sw_p^p(P_n',Q)$, where  $P_n'= \frac{1}{n}( \delta_{X_1'} + \sum_{i=2}^n \delta_{X_i})$, with $X_1'$ drawn from $P$, independent of $X_1,\ldots,X_n$. To simplify the notation, let $X_j'=X_j$ for $j\geq 2$.\\

For each $\theta \in \mathbb S^{d-1}$, consider the following construction. First, define the random map $\tau'_\theta:\{1,\ldots,n\} \rightarrow \{1,\ldots,n\}$ that maps each index $i$ to the position of $\langle\theta,X_i'\rangle$ in the ordered statistic associated with the variables $\{\langle\theta,X_i'\rangle\}_{i=1}^n$. If there are $k$ equal values in the ordered statistic, then $\tau'_\theta$ is randomly selected as one of the $k!$ possible permutations. With this construction, $   \tau'_\theta(1),\ldots , \tau'_\theta(n)$ are identically distributed, 
and given that  $X_1$ is independent of $X_1'\ldots,X_n'$, 
\begin{equation}\label{eq:equally_distr}
    (X_1,\tau'_\theta(1))\overset{d}{=}\ldots \overset{d}{=}(X_1,\tau'_\theta(n)) \ .
\end{equation}
Denote by $F_{\theta,n},F_{\theta,n}',G_\theta$ the distribution functions  of $P_n^\theta,P_n'^\theta$ and $Q^\theta$, respectively. Given a random variable $U$ with uniform distribution in $(0,1)$, define the positive measure 
$    Q_j^\theta(A) = \mathbb P\bigl( G_\theta^{-1}(U)\in A, U \in (\frac{j-1}{n},\frac{j}{n}] \bigr)
$  for each $j=1,\ldots,n$.
Note that $Q^\theta = \sum_{j=1}^n Q_j^\theta$, and $Q_j^\theta(\mathbb R)=1/n$. Moreover, by the closed-form expression for $W_p$ in terms of quantile functions,   $\pi_n'^\theta=\mathcal L(F_{\theta,n}'^{-1}(U),G_\theta^{-1}(U))$ is an optimal coupling between $P_n'^\theta$ and $Q^\theta$, and can be decomposed as $ \pi_n'^\theta =\sum_{i=1}^n \delta_{\langle \theta, X_i' \rangle} \times Q_{\tau'_\theta(i)}^\theta $. Similarly, $\pi_n^\theta = \sum_{i=1}^n \delta_{\langle\theta,X_i\rangle} \times Q^\theta_{\tau_\theta'(i)}$ is a valid coupling between $P_n^\theta$ and $Q^\theta$.
Therefore, since $X_i=X_i'$ for $i\neq 1$,
\begin{align*}
\w_p^p(P_n^\theta,Q^{\theta}) -\w_p^p(P_n'^\theta,Q^{\theta}) \leq  \int_{\mathbb R} |\langle \theta,X_1\rangle-t|^p -|\langle \theta,X_1'\rangle-t|^p dQ^\theta_{\tau'_\theta(1)}(t)
\end{align*} 
which, using the convexity of $x \mapsto |\langle \theta,x\rangle-t|^p$ for every $\theta\in\mathbb S^{d-1}$, ensures that
\begin{align} \label{eq:primal_inequality}
    Z_n-Z_n'  &\leq p\| X_1 -X_1'\| \int_{\mathbb S^{d-1}}  \int_{\mathbb R}  |\langle \theta,X_1\rangle-t|^{p-1}  dQ^\theta_{\tau'_\theta(1)}(y)  d\sigma(\theta) \ .
\end{align}
Applying Hölder's inequality, with $f=1,g=|\langle \theta,x\rangle-t|^{p-1}, q_1=\frac{2p}{p+1}$, $q_2=\frac{2p}{p-1}$, 
\begin{align}\label{eq:des1}
    \int_{\mathbb R}  |\langle \theta,X_1\rangle-t|^{p-1}  dQ^\theta_{\tau'_\theta(1)}(t)  &\leq \frac{1}{n^\frac{p+1}{2p}}\left( \int_{\mathbb R} |\langle \theta,X_1\rangle-t|^{2p} dQ^\theta_{\tau_\theta'(1)}(t) \right)^{\frac{p-1}{2p}} \ .
\end{align}
By exchangeability \eqref{eq:equally_distr}, we know that for every $\theta \in \mathbb S^{d-1}$, 
\begin{equation}
    \label{eq:des2}
   \hspace{-0.8cm} \mathbb E \left( \int_{\mathbb R} |\langle \theta,X_1\rangle-t|^{2p} dQ^\theta_{\tau_\theta'(1)}(t)\right)= \frac{1}{n} \mathbb E \left( \sum_{i=1}^n\int_{\mathbb R} |\langle \theta,X_1\rangle-t|^{2p}dQ_{\tau_\theta'(i)}^\theta(t) \right) 
     =  \frac{1}{n}  \int_{\mathbb R^d} |\langle \theta,X_1\rangle-t|^{2p} dQ^\theta(t) \ .
\end{equation}
From \eqref{eq:primal_inequality} and a further application of Hölder's inequality with $q_1=p$, $q_2=\frac{p}{p-1}$,
\begin{align}\label{eq:Hölder_proof1}
    \mathbb E &\Bigl( (Z_n-Z_n')_+^2 \Bigr)   \leq p^2 \mathbb E \Bigl(   \| X_1 -X_1'\|^{2p}\Bigr)^\frac{1}{p}\mathbb E \Biggl( \Bigl( \int_{\mathbb S^{d-1}}  \int_{\mathbb R}  |\langle \theta,X_1\rangle-t|^{p-1}  dQ_{\tau'_\theta(1)}(y)  d\sigma(\theta) \Bigr)^\frac{2p}{p-1} \Biggr)^\frac{p-1}{p} \ .
\end{align}
Applying Jensen's inequality, \eqref{eq:des1} and \eqref{eq:des2}, we obtain the bound
\begin{align*}
    &\mathbb E \Biggl( \Bigl( \int_{\mathbb S^{d-1}}  \int_{\mathbb R}   |\langle \theta,X_1\rangle-t|^{p-1}  dQ_{\tau'_\theta(1)}(y)  d\sigma(\theta) \Bigr)^\frac{2p}{p-1}\Biggr) =\\
    & \leq \frac{1}{n^\frac{2p}{p-1}} \mathbb E \Biggl(  \int_{\mathbb S^{d-1}} \int_{\mathbb R} |\langle \theta,X_1\rangle-t|^{2p} dQ^\theta(t) d\sigma(\theta) \Biggr) \leq \frac{1}{n^\frac{2p}{p-1}} \mathbb E \Bigl(  \|X_1-Y\|^{2p} \Bigr) \ ,
\end{align*}
where $Y$ represents a random variable independent of $X$, with $\mathcal L(Y)=Q$. Combined with \eqref{eq:Hölder_proof1}, we deduce that $\mathbb E \bigl( (Z_n-Z_n')_+^2\bigr) \leq \frac{1}{n^2}K(P,Q)$, which allows us to conclude (a) using the Efron-Stein inequality.  \\

\noindent
\textit{Two-sample setting.} Let $Z= \sw_p^p(P_n,Q_m)$, where $Q_m$ is the empirical distribution associated with $m$ i.i.d. random variables $Y_1,\ldots,Y_m$ drawn form $Q$, independent of $X_1,\ldots,X_n$. $Z$ is symmetric with respect to its first $n$ variables and its last $m$ variables. Thus, applying Efron-Stein,
\begin{equation*}
    \operatorname{Var}(\sw_p^p(P_n,Q_m)) \leq n \mathbb E\bigl( (Z_{n,m}-Z_{n,m}')_+^2\bigr) + m \mathbb E\bigl( (Z_{n,m}-Z_{n,m}'')_+^2\bigr) \ ,
\end{equation*}
where $Z_{n,m}'=\sw_p^p(P_n',Q_m)$, $Z_{n,m}''=\sw_p^p(P_n,Q_m')$, $P_n'$ is defined as in the one-sample case, and $Q_m'$ is the empirical measure associated with $Y_1',Y_2\ldots,Y_m$, with $Y_1'$ an independent copy. Denote by $\mathbb E_X$ the expectation taken only with respect to $X_1,\ldots,X_n$. Applying the one-sample reasoning to the discrete  $Q=Q_m$ ensures that 
\begin{align*}
    \mathbb E_X \bigl( (Z_{n,m}-Z_{n,m}')_+^2\bigr)  \leq \frac{p^2}{n^2} \mathbb E_X \Bigl(   \| X_1 -X_1'\|^{2p}\Bigr)^\frac{1}{p}  \mathbb E_X \Bigl(  \frac{1}{m}\sum_{j=1}^m\|X_1- Y_j\|^{2p} \Bigr) ^\frac{p-1}{p} \ .
\end{align*}
Taking the expectation w.r.t.  $Y_1,\ldots,Y_n$ and applying Jensen's inequality yields $\mathbb E \bigl( (Z_{n,m}-Z_{n,m}')_+^2\bigr) \leq \frac{1}{n^2}K_p(P,Q)$. Similarly, exchanging the roles of $P$ and $Q$, $\mathbb E\bigl( (Z_{n,m}-Z_{n,m}'')_+^2\bigr)\leq \frac{1}{m^2}K_p(Q,P)$, which allows us to conclude (b). 
\qed

\subsection{Proof of Proposition \ref{teor:var_dual}}

Let $P_n,P_n'$ and  $Z_n,Z_n'$ be as in the proof of \Cref{teor:primal}, and define $R_n'$ as $R_n$, replacing $P_n$ with $P_n'$. Let $\phi_{n}^\theta$ be any $c$-concave optimal transport potential from $P_n^\theta$ to $Q^\theta$, and let  $\phi^\theta=\phi^\theta_{x_0}$.
The first step in the proof consists of showing that $n(R_n-R_n')_+\rightarrow 0$ $\mathbb P$-a.s. To that end, consider the set of $\mathbb P$-probability one $\Omega_0\subset \Omega$ where all the properties of Proposition \ref{prop:regularidad_potenciales} hold, and where $X_1' \in \operatorname{int}(\supp P)$. For any $\omega\in \Omega_0$, given $\{X_n^{\omega}\}_{n=1}^\infty$ and $X_1'^\omega$, the goal is to demonstrate that $n(R_n^\omega-R_n'^\omega)_+\rightarrow 0$. For simplicity, we omit the superscript $\omega$, noting that, within the considered set $\Omega_0$, all assertions of Proposition \ref{prop:regularidad_potenciales} are satisfied. In particular,  $\phi_{n,X_1}^\theta$  and $\phi^\theta_{X_1}$, defined in \Cref{prop:regularidad_potenciales}, respectively, are well-defined $c$-concave optimal transport potentials satisfying 
$\phi_{n,X_1}^\theta(\langle \theta,X_1\rangle)=\phi_{X_1}^\theta(\langle \theta,X_1\rangle)=0$  for every $\theta\in\mathbb S^{d-1}$ and $n\in\mathbb N  .$
Reordering terms in \eqref{eq:linearization} and integrating over $\mathbb S^{d-1}$, it follows that 
\begin{align*}
 R_n-R_n' \leq&\int_{\mathbb S^{d-1}}  \int_{\mathbb R} (\phi^\theta_{n,X_1}(s)-\phi^\theta_{X_1}(s))d(P_n^\theta- P_n'^\theta)(s)d\sigma(\theta) \\
= & \ \frac{1}{n}  \int_{\mathbb S^{d-1}}  \Bigl( \phi_{n,X_1}^\theta( \langle \theta, X_1 \rangle) - \phi^\theta_{X_1}( \langle \theta, X_1 \rangle )
 \Bigr)d\sigma(\theta) -  \frac{1}{n}  \int_{\mathbb S^{d-1}}  \Bigl( \phi_{n,X_1}^\theta( \langle \theta, X_1' \rangle) - \phi^\theta_{X_1}( \langle \theta, X_1' \rangle )
 \Bigr)d\sigma(\theta) 
\end{align*}
Therefore, to show that $n(R_n-R_n')_+\rightarrow 0$, it suffices to prove  that
\begin{equation}\label{eq:con_integrales}
    \int_{\mathbb S^{d-1}}  \phi_{n,X_1}^\theta( \langle \theta, x \rangle )d\sigma(\theta) \longrightarrow \int_{\mathbb S^{d-1}}  \phi_{X_1}^\theta( \langle \theta, x \rangle )d\sigma(\theta)  \quad   \forall \ x \in \operatorname{int}(\supp P) \ .
\end{equation}
By Proposition \ref{prop:regularidad_potenciales}, for every $x \in \operatorname{int}(\operatorname{supp}(P))$ and   $\theta\in\mathbb S^{d-1}$, $\phi_{n,X_1}^\theta( \langle \theta, x\rangle) \rightarrow \phi_{X_1}^\theta( \langle \theta, x\rangle)$. Moreover, there exist $n_0\in\mathbb N$ and $M_{x}>0$ such that $|\phi_{n,X_1}^\theta( \langle \theta, x\rangle)| \leq M_{x}$ for every  $n\in \mathbb N$ and $\theta \in \mathbb S^{d-1}$. Therefore, \eqref{eq:con_integrales} follows from the  dominated convergence theorem, concluding that  $n(R_n-R_n')_+\rightarrow 0$ pointwise in $\Omega_0$. 
Our goal now is to check that $\mathbb E (n^2(R_n-R_n')^2_+) \rightarrow 0$. Given  that $n(R_n-R_n')_+\rightarrow 0 \ \mathbb P$-a.s., it follows that $n^2(R_n-R_n')^2_+\rightarrow 0 \ \mathbb P$-a.s. Consequently, it suffices to prove uniform integrability of $n^2(R_n-R_n')^2_+$. Recall that
\begin{align}\label{eq:partes_unif_int}
    n(R_n-R_n') =  n(Z_n-Z_n') - \int_{\mathbb S^{d-1}} \phi^\theta_{x_0}(\langle \theta,X_1\rangle) d\sigma(\theta)  +  \int_{\mathbb S^{d-1}} \phi^\theta_{x_0}(\langle \theta,X_1'\rangle)\bigr) d\sigma(\theta)
\end{align}
The last two terms do not depend on $n$, and their second-order moment is bounded by $\|\phi^\theta_{x_0}\|_{L^2(\sigma\times P)}^2$, which is finite by assumption \eqref{eq:int_potentials_limit_x0}. Thus, it suffices to show uniform integrability of $n^2(Z_n-Z_n')^2_+$, for which it is enough to prove uniform integrability of the square of $n$ times the right-hand side in \eqref{eq:primal_inequality}. To that end, we provide a uniform bound on the $1+\gamma$ moment. With the moment assumptions of \Cref{teor:var_dual}, this can be done exactly with the same arguments as in the proof of \Cref{sec:primal}, adapting the powers in the H\"older's inequality in the same style of the proof of Theorem 4.6 in \cite{delBarrio2024CLTgeneral}.
Therefore,  $\mathbb E (n^2(R_n-R_n')^2_+) \rightarrow 0$, and Efron-Stein inequality allows us to conclude.\\

\noindent
\textit{Two-sample setting.} Let $P_n,P_n',Q_m,Q_m'$ be defined as in the proof of Proposition \ref{teor:var_dual}. Define $R_{n,m}'$ by replacing $P_n$ with $P_n'$ in the definition of $R_{n,m}$. Similarly, define $R_{n,m}''$ by replacing $Q_m$ with $Q_m'$. Applying Efron-Stein inequality,
\begin{equation*}
   \frac{nm}{n+m} \operatorname{Var}(R_{n,m}) \leq n^2 \mathbb E\bigl( (R_{n,m}-R_{n,m}')_+^2\bigr) + m^2 \mathbb E\bigl( (R_{n,m}-R_{n,m}'')_+^2\bigr) \ .
\end{equation*}
The convergence $ n^2 \mathbb E\bigl( (R_{n,m}-R_{n,m}')_+^2\bigr) \rightarrow 0$, follows from the same argument as in the one-sample setting, considering now any $c$-concave  optimal transport potential  $\phi_{n,m}^\theta$ from $P_n^\theta$ to $Q_m^\theta$, and defining $\phi_{n,m,X_1}^\theta$  as in \Cref{prop:regularidad_potenciales} and $\phi_{X_1}^\theta$ as before. Property (6) in Proposition \ref{prop:regularidad_potenciales} and \Cref{prop:integrability_potentials} enable us to extend the previous arguments for the almost sure convergence, whereas the uniform integrability follows from $\|\phi^\theta_{x_0}\|_{L^2(\sigma\times P)}^2<\infty$, the previous bounds, and the same argument as in the proof of Proposition \ref{teor:primal} for the two-sample setting.
Similarly, one can prove $m^2 \mathbb E(R_{n,m}-R_{n,m}'')_+^2 \rightarrow 0 $. Again, we need to consider a suitable choice of $c$-concave optimal transport potentials, but now fixed at $Y_1$. Given any $c$-concave optimal transport potential $\psi_{n,m}^\theta$ from $Q_m^\theta$ to $P_n^\theta$ (for instance,  $\psi_{n,m}^\theta= (\phi_{n,m}^\theta)^c$), define
$\psi_{n,m,Y_1}^\theta = \psi_{n,m}^\theta - \psi_{n,m}^\theta(\langle \theta, Y_1 \rangle) $ and $  \psi_{Y_1}^\theta= (\phi_{x_0}^\theta)^c - (\phi_{x_0}^\theta)^c (\langle \theta, Y_1 \rangle)$.
Taking into account the symmetric role of $P$ and $Q$,  property (6) in  Proposition \ref{prop:regularidad_potenciales} is verified, exchanging the roles of $P$ and $Q$. The same argument as before allows us to conclude, now using that, by \eqref{eq:int_potentials_limit_x0}, $\|(\phi^\theta_{x_0})^c\|_{L^2(\sigma\times Q)}^2<\infty$.

\qed

\subsection{Proof of Theorem \ref{teor:CLT}}
\noindent
With the slight abuse of notation presented in \Cref{sect:CLT_general}, if $\bar Y_n= \frac{1}{n}\sum_{i=1}^nY_i$, 
 \Cref{teor:var_dual} ensures that
    \begin{equation*}
         \mathbb E \left(\Bigl( \sqrt n \bigl(\sw_p^p(P_n,Q)-\mathbb E \sw_p^p(P_n,Q)\bigr) - \sqrt{n}\bigl(\bar Y_n-\mathbb EY\bigr)\Bigr)^2\right) \rightarrow 0 \ ,
    \end{equation*}  
which implies convergence in the Wasserstein distance $W_2$.  Since the convergence $\sqrt{n}(\bar Y_n-\mathbb EY)\rightsquigarrow N(0,v_{P,Q}^2)$ was obtained from the CLT, we have convergence of second-order moments, and therefore convergence in $\w_2$ (see Lemma 8.3 in \cite{bickel1981someAsymptotic}). Applying the triangle inequality, we conclude
$   W_2 \bigl( \mathcal L\bigl(\sqrt n (\sw_p^p(P_n,Q)-\mathbb E (\sw_p^p(P_n,Q))\bigr),N(0,v_{P,Q}^2) \bigr) \rightarrow 0 \ .
$ \\

\noindent
\textit{Two-sample setting.} The exact same reasoning applies. \qed

\section{Proofs of Section \ref{sect:bias}}\label{appendixB_bueno}

\subsection{Additional lemmas}
\noindent
Throughout this section, we repeatedly rely on the following elementary inequality.  Given $x,y\in \mathbb R^d$ and $q > 0$, 
\begin{equation}\label{eq:elemental_ineq_p}
    \| x + y \|^q \leq C_q \bigl( \|x\|^q + \|y\|^q)\, \quad \quad \textnormal{where } \left\{\begin{array}{lc}
       C_q=1  &  0<q\leq 1\\
       C_q = 2^{q-1}  &  q\geq 1
    \end{array}\right. \ .
\end{equation}
First, we present a well-known fact about the convergence in mean of the Wasserstein distance (see Theorem 2.14 in \cite{bobkov2019onedimensional}), for which we include a shorter alternative proof.

\begin{lemma}\label{lemma:convergence_EW_p}
    Let $p\geq 1$. For any $P\in\mathcal P_p(\mathbb R)$, $\mathbb E\bigl(\w_p^p(P_n,P)\bigl)\rightarrow 0$.
\end{lemma}

\begin{proof}
The proof is based on studying  $|F_n^{-1}(\omega,t)-F^{-1}(t)|^p$ as a random variable in the product space $\Omega\times (0,1)$, with the product probability $\mathbb P\times \ell_1$, where $\ell_1$ denotes the Lebesgue measure. The Glivenko-Cantelli theorem ensures that there exists $\Omega_0\subset\Omega$ with $\mathbb P(\Omega_0)=1$ such that $\sup_{x\in\mathbb R} |F_n(\omega,x)-F(x)| \rightarrow 0$, for every  $\omega\in\Omega_0$. Thus, for any of such $\omega\in \Omega_0$, $ P_n^\omega\rightsquigarrow P$. If $C_{F^{-1}}$ denotes the set of continuity points of $F^{-1}$, which verifies $\ell_1(C_{F^{-1}})=1$,  by the Skorohod elementary theorem (e.g. see 
pp. 9-10 in \cite{shorack2009empirical}),  
\begin{equation}\label{eq:elementary_skorohod}
   |F_n^{-1}(\omega,t)-F^{-1}(t)|^p\rightarrow 0 \quad \forall (\omega,t)\in\Omega_0 \times C_{F^{-1}} \ .
\end{equation}
Thus, we have almost sure convergence in $\Omega\times(0,1)$, and the only thing left to prove is uniform integrability of the sequence. From the inequality \eqref{eq:elemental_ineq_p}, it suffices to prove uniform integrability of $|F_n^{-1}(\omega,t)|^p$, which follows from
\begin{align}\label{eq.distr_Fn^-1}
  \mathbb P\times\ell_1\Bigl( \{ (\omega,t):F_n^{-1}(\omega,t)\leq x\}\Bigr)  &= \mathbb P\times\ell_1\Bigl( \{ (\omega,t):t\leq F_n(\omega,x)\}\Bigr)  = \int_\Omega F_n(\omega,x) d\mathbb P(\omega) = F(x) \ . 
\end{align}
Therefore, as a random variable in $(0,1)\times\Omega$, $\mathcal L(F_n^{-1})=P\ \forall\  n \in \mathbb N$, which allows us to conclude.
\end{proof}

Given $p> 1$, the function $h_p(x)=|x|^p$ is differentiable, with derivative  $h_p'(x)= p\ \sg(x)|x|^{p-1}$. The following lemma bounds the expected values of integrals involving the difference between $h'_p$ evaluated at the empirical quantile function difference and its population counterpart. Although formulated for general $p>1$, improved assumptions could be given for the simpler case $p=2$.
Since the proof of Theorem \ref{teor:bias_improved_csorgo_sliced} for $p=2$ does not rely on this lemma, we do not pursue this direction further.
\begin{lemma}\label{lemma:convergence_derivative}
    Let $p>1$, $\beta\geq 1$ and $P,Q\in\mathcal P_{\beta(p-1)}(\mathbb R)$, with distribution functions $F$ and $G$, respectively. Then
\begin{align}
        &\mathbb E \Bigl( \int_0^1 \big| h'_p(F_n^{-1}(t)-G^{-1}(t)) - h'_p(F^{-1}(t)-G^{-1}(t)) \big|^\beta dt \Bigr) \rightarrow 0 \ , \label{eq:bound_derivative_hp} \\
        &\mathbb E \Bigl( \int_0^1 \big| h'_p(F_n^{-1}(t)-G_m^{-1}(t)) - h'_p(F^{-1}(t)-G^{-1}(t)) \big|^\beta dt \Bigr) \rightarrow 0 \ . \label{eq:bound_derivative_hp_2}
\end{align}
\end{lemma}

\begin{proof} 
From the inequality \eqref{eq:elemental_ineq_p}, it follows that 
\begin{align}
     &\frac{1}{p C_\beta^\beta}\Big| h'_p(F_n^{-1}(t)-G^{-1}(t)) - h'_p(F^{-1}(t)-G^{-1}(t)) \Big|^\beta \notag \\
     & \leq  \big| \sg(F_n^{-1}(t)-G^{-1}(t))\bigl( |F_n^{-1}(t)-G^{-1}(t)|^{p-1} - |F^{-1}(t)-G^{-1}(t)|^{p-1} \bigr)\Big|^\beta  \label{eq:lemmaB2_1}\\
     &+\Big|  \bigl( \sg(F_n^{-1}(t)-G^{-1}(t)) - \sg(F^{-1}(t)-G^{-1}(t)) \bigr)  |F^{-1}(t)-G^{-1}(t)|^{p-1}\Big|^\beta \ . \label{eq:lemmaB2_2}
\end{align}
With the same reasoning as in \eqref{eq:elementary_skorohod}, $F_n^{-1}(t) \rightarrow F^{-1}(t)$  for every  $t\in C_{F^{-1}}$ in a set of $\mathbb P$-probability one. Therefore, \eqref{eq:lemmaB2_2} converges to zero almost surely as a random variable in $\Omega\times(0,1)$, and it is bounded by the integrable function $2^\beta |F^{-1}(t)-G^{-1}(t)|^{\beta(p-1)}$. Thus, applying dominated convergence theorem we can conclude that the expected value of the integral of \eqref{eq:lemmaB2_2} converges to zero.  
Now, the expected value of the integral of \eqref{eq:lemmaB2_1} is bounded by 
\begin{equation}
    \mathbb E \Bigl( \int_0^1\Bigl| |F_n^{-1}(t)-G^{-1}(t)|^{p-1} - |F^{-1}(t)-G^{-1}(t)|^{p-1} \Big|^\beta dt \Bigl) \label{eq:lemmaB2_3}
\end{equation}
We distinguish two different cases: \\

\noindent
$\bullet$ If $1<p\leq 2$, \eqref{eq:elemental_ineq_p} implies that $\big||x|^{p-1} - |y|^{p-1}\big| \leq |x-y|^{p-1}$. Thus, \eqref{eq:lemmaB2_3} is bounded by 
$\mathbb E \bigl( \w_{\beta(p-1)}^{\beta(p-1)}(P_n,P)\bigr)$,
which converges to zero by Lemma \ref{lemma:convergence_EW_p}.\\

\noindent
$\bullet$ If $p> 2$, then $h_{p-1}(\cdot)= |\cdot|^{p-1}$ is convex and differentiable. Thus, we know that for every $x,y\in\mathbb R$
\begin{align}\label{eq:convex_basic_ineq}
    & h_{p-1}'(y)(x-y) \leq h_{p-1}(x) -h_{p-1}(y) \leq h_{p-1}'(x)(x-y) \ .
\end{align}
Using this inequality,  we can bound \eqref{eq:lemmaB2_3}, up to a constant factor, by 
\begin{align*}
     & \mathbb E \Bigl(\int_0^1   | F_n^{-1}(t) - G^{-1}(t)|^{\beta(p-2)}| F_n^{-1}(t) - F^{-1}(t) |^\beta dt \Bigr)  +  \mathbb E \Bigl(\int_0^1   \bigl| F^{-1}(t) - G^{-1}(t)|^{\beta(p-2)}| F_n^{-1}(t) - F^{-1}(t) |^\beta dt \Bigr) \ .
\end{align*}
 We prove that the first term converges to zero, for the second term follows similarly. Applying Hölder's inequality with $q_1=\frac{p-1}{p-2}$, $q_2=p-1$, the first term is bounded by 
\begin{align*}
\mathbb E \Bigl(\int_0^1   | F_n^{-1}(t) - G^{-1}(t)|^{\beta(p-1)}| dt \Bigr)^{\frac{p-2}{p-1}}      \mathbb E \Bigl(\int_0^1   | F_n^{-1}(t) - F^{-1}(t)|^{\beta(p-1)} dt \Bigr)^{\frac{1}{p-1}} \ .
\end{align*}
By the moment assumptions, the first term is bounded, whereas the second term converges to zero by Lemma \ref{lemma:convergence_EW_p}. \\

\noindent
\textit{Two-sample setting:} The proof follows from the same steps, and a further application of the inequality \eqref{eq:elemental_ineq_p}.
\end{proof}

\subsection{Proof of \Cref{lemma:csorgo}}

Let $X_1,\ldots,X_n$ be an i.i.d. sample from $P$, and assume that $X_i=F^{-1}(Z_i)$, where $Z_1,\ldots,Z_n$ are i.i.d. $U(0,1)$.  If $U_n$ denotes the cumulative distribution function of $Z_1,\ldots,Z_n$, the weighted general quantile process and the uniform quantile process are defined for each $t\in(0,1)$ as
\begin{align*}
    \rho_n(t) := f(F^{-1}(t)) \sqrt{n}(F_n^{-1}(t)-F^{-1}(t)) \ ,\hspace{0.8cm}    u_n(t) := \sqrt{n}(U_n^{-1}(t)-t) \ .
\end{align*}
Denoting $r(t)= \frac{h(t)}{f(F^{-1}(t))}$, our goal is to show
$\int_{1/n}^{1-1/n} \rho_n(t) r(t) dt \rightsquigarrow 
    N(0,v^2).$ 
The proof is based on two steps. First, we show that the convergence holds if we replace the weighted general quantile process by the uniform quantile process. Then, we demonstrate that the induced error converges in probability to zero. For the first part, following \cite{delBarrio2005asyptotics}, if $\{\xi_n\}_{n=1}^\infty$ denotes a sequence of i.i.d. exponential variables with mean 1, then 
\begin{equation}\label{eq:exponentials}
  \int_{1/n}^{1-1/n} u_n(t) r(t) dt  \overset{d}{=} \frac{n}{S_{n+1}} \sum_{i=1}^{n+1} C_{n,i}\xi_i = \frac{n}{S_{n+1}} \Bigl(\sum_{i=1}^{n+1} C_{n,i}(\xi_i-1) + \sum_{i=1}^{n+1} C_{n,i}\Bigr) \ ,   
\end{equation}
where $C_{n,i}=n^{-1/2}\int_{1/n}^{1-1/n} r(t)a_{n,i}(t)dt$ and $a_{n,i}(t)= (1-t)I(t>\frac{i-1}{n}) - tI(t\leq\frac{i-1}{n})$. Let $\Tilde m_n(t) = \sum_{i=1}^{n+1} a_{n,i}(t)$, $\Tilde K_n(s,t) = \sum_{i=1}^{n+1} a_{n,i}(s)a_{n,i}(t)$. Note that Assumptions (i), (ii) and (iii) together with H\"older's inequality ensure that 
$\frac{1}{\sqrt{n}}\int_{1/n}^{1-1/n} |r(t)|dt\to 0$
and 
$ \int_0^1\int_0^1  |r(s)r(t)|(s\wedge t-st) dt ds<\infty$. 
Then, using Lemma 2.1 in \cite{delBarrio2005asyptotics}, it follows that 
\begin{itemize}
    \item $ |\sum_{i=1}^{n+1} C_{n,i}| = |\frac{1}{\sqrt{n}}\int_{1/n}^{1-1/n} r(t)\Tilde m_n(t)dt|\leq \frac{1}{\sqrt{n}}\int_{1/n}^{1-1/n} |r(t) |dt\rightarrow 0 \ ,$
    \item $\max_i |C_{n,i}|\leq \frac{1}{\sqrt{n}}\int_{1/n}^{1-1/n} |r(t) |dt\rightarrow 0  \ ,$
    \item $\sum_{i=1}^{n+1} C_{n,i}^2 =\int_{1/n}^{1-1/n}\int_{1/n}^{1-1/n}r(s)r(t)\frac{\Tilde K_n(s,t)}{n}dsdt \rightarrow v^2 .$
\end{itemize}
Then, applying Lemma 3.7 in \cite{delBarrio2005asyptotics} and the strong law of large numbers we obtain $\int_{1/n}^{1-1/n} u_n(t) r(t) dt \rightsquigarrow 
    N(0,v^2).$ 
Therefore, our goal now is to show that $\int_{1/n}^{1-1/n} \rho_n(t) r(t) dt - \int_{1/n}^{1-1/n} u_n(t) r(t) dt \rightarrow_p 0$, for which it suffices to prove that
\begin{itemize}
    \item For every $0<\epsilon<1/2$, $\int_{\epsilon}^{1-\epsilon} \rho_n(t) r(t) dt - \int_{\epsilon}^{1-\epsilon} u_n(t) r(t) dt \rightarrow_p 0  $.
    \item $\lim_{\epsilon\downarrow 0} \limsup_{n\rightarrow\infty }\mathbb P( |\int_{1/n}^{\epsilon} u_n(t) r(t) dt |>\delta)= \lim_{\epsilon\downarrow 0} \limsup_{n\rightarrow\infty } \mathbb P(|\int_{1/n}^{\epsilon} \rho_n(t) r(t) dt |>\delta) =  0$ for every $\delta>0$, and similarly for $\int_{1-\epsilon}^{1-1/n}$.
\end{itemize}
The argument follows ideas similar to those in the proof of Theorem 2.1 in \cite{csorgo1990distributions},
but without relying on the strong approximations of \cite{csorgo1978strongApproximations}.
First, applying the mean value theorem, for some $\eta(t)\in [U_n^{-1}(t)\wedge t, U_n^{-1}(t)\vee t]$, we have
\begin{equation*}
   \big| \int_{\epsilon}^{1-\epsilon} (\rho_n(t) - u_n(t)) r(t) dt \big| \leq \sup_{\epsilon\leq t \leq 1-\epsilon} \big| \frac{f(F^{-1}(t))}{f(F^{-1}(\eta(t)))}-1 \bigl| \big| \int_{\epsilon}^{1-\epsilon}  u_n(t) r(t) dt \big| \rightarrow 0 \ .
\end{equation*}
This convergence follows from the convergence in probability to zero of the first term (as in \cite{csorgo1990distributions}) together with the weak convergence of the second term, which can be shown using same reasoning as in the first part of the proof. Similarly, it follows that $\int_{1/n}^{\epsilon} u_n(t) r(t) dt \rightsquigarrow 
    N(0,v_\epsilon^2)$ where  $v_\epsilon^2 := \int_0^\epsilon\int_0^\epsilon  \frac{h(t)h(s)}{f(F^{-1}(t))f(F^{-1}(s))}(s\wedge t-st) dtds$, which satisfies $v_\epsilon^2\to 0$ by the dominated convergence theorem. This shows that, for every $\delta>0$,  $\limsup_{n\rightarrow\infty }\mathbb P( |\int_{1/n}^{\epsilon} u_n(t) r(t) dt |>\delta)= 2\Phi(-\delta/v_\epsilon) \rightarrow 0$ as $\epsilon \rightarrow 0$, where $\Phi$ denotes the distribution function of a standard normal, which gives the first convergence in the second property. To prove the second convergence, reasoning as in \cite{csorgo1990distributions}, it suffices to prove $\lim_{\epsilon\downarrow 0} \limsup_{n\rightarrow\infty }\mathbb P( |\int_{1/n}^{\epsilon} \frac{1}{f(F^{-1}(t/\lambda))}u_n(t) h(t) dt |>\delta)=0$ for every $\lambda>1$ if $f\circ F^{-1}$ is non-increasing around 0, or for every $\lambda<1$ if it is non-decreasing. This follows if we show that the sequence of integrals weakly converge to a centered Gaussian for every $0<\epsilon<\lambda/2\wedge1/2$, and the  asymptotic variances converge to zero as $\epsilon \to 0$. Again, this can be shown as above, now using that if we denote $\tilde r_\epsilon(t)= I(0\leq t \leq \epsilon)h(t)/f(F^{-1}(t/\lambda))$, applying H\"older's inequality and using the change of variable formula, 
 $n^{-1/2}\int_{1/n}^{1-1/n} |\tilde r_\epsilon(t)|dt  \to 0$ as $n\rightarrow \infty$, and using also that $s\wedge t-\lambda st\leq 2 (s\wedge t-st)$ for every $t\in (0,1/2)$ and $\lambda>0$,
 $ \int_0^1\int_0^1  |r_\epsilon(s)r_\epsilon(t)|(s\wedge t-st) dt ds<\infty$ and converges to zero as $\epsilon\to 0$.\\

Finally, for the second part of the proposition, we only need to prove that adding $\int_0^{1/n}$ (or symmetrically $\int_{1-1/n}^1$) does not change the asymptotic behavior.  This follows from Markov's inequality and,  if $J_\alpha(P)<\infty$, 
\begin{align*}
    \mathbb E \Bigl( \Big|\int_0^{\frac{1}{n}} h(t)\sqrt{n}(F_n^{-1}(t)-F^{-1}(t))dt  \Big| \Bigr)& \leq \Bigl( \int_0^{\frac{1}{n}} |h(t)|^\beta dt \Bigr)^{1/\beta}  \Bigl( n^{\alpha/2} \mathbb E (W_\alpha^{\alpha}(P_n,P))\Bigr)^{1/\alpha}  \rightarrow 0 \ ,
\end{align*}
where the convergence follows from Assumption (iii) and  Proposition 5.3 in \cite{bobkov2019onedimensional}. If $P\in\mathcal P_{2p}(\mathbb R)$ and $\int_0^1 |h(t)|^{2p/(p-1)} dt<\infty$, then  $\|h F^{-1}\|_{L^2(0,1)}<\infty$ and we can conclude as in Lemma 2.2 in \cite{alvarez2011uniqueness}.

\subsection{Proof of \Cref{teor:bias_improved_csorgo_sliced}}

\vspace{0.1cm}
\noindent
\textbf{Case $\mathbf{d=1}$ \textbf{:}\ }
For simplicity, we prove first Theorem \ref{teor:bias_improved_csorgo_sliced} in the one-dimensional setting. For this, denote by $f$ and $F$ the density and distribution function of $P$, which play the role of $f_\theta,F_\theta$ in the statement of Theorem \ref{teor:bias_improved_csorgo_sliced}, and note that $\operatorname{J}_{\alpha}(P) = \operatorname{SJ}_{\alpha}(P)<\infty$.
Let $F_n, G, G_m$ be the distribution function of $P_n, Q$ and $Q_m$, respectively.
For any $p>1$, with $p\neq 2$, if we denote $T_n  =  \sqrt{n} \bigl(\w_p^p(P_n,Q) - \w_p^p(P,Q)\bigr) $, then  by \eqref{eq:convex_basic_ineq} and the quantile expression, if we define  
\begin{align*}
     &L_n = \int_0^1 h'_p\bigl( F^{-1}(t) - G^{-1}(t)\bigr)\sqrt{n}\bigl( F_n^{-1}(t) - F^{-1}(t) \bigr) dt 
    \ , \quad \\
     &U_n=  \int_0^1 h
     '_p\bigl( F_n^{-1}(t) - G^{-1}(t)\bigr)\sqrt{n}\bigl( F_n^{-1}(t) - F^{-1}(t) \bigr) dt  \ ,
\end{align*}
then $L_n\leq T_n\leq U_n$. Under the given assumptions, it follows from H\"older's inequality together with \Cref{lemma:convergence_derivative}  and Proposition 5.3 in \cite{bobkov2019onedimensional} that $U_n-L_n \rightarrow_p 0$. Using \Cref{lemma:csorgo} with $h(t) = h'_p(F^{-1}(t)-G^{-1}(t))$, it is immediate to show that $L_n$ (and thus, $T_n$) weakly converges to a centered Gaussian distribution. Finally, under the given assumptions, \Cref{teor:primal}  ensures that  the sequence $\sqrt{n}(\w_p^p(P_n,Q)  - \mathbb E(\w_p^p(P_n,Q)))$ is uniformly integrable. Then, a tightness argument implies boundedness of  $\sqrt{n}(\mathbb E(\w_p^p(P_n,Q))-\w_p^p(P,Q) )$ and hence, uniform integrability of
\begin{align}
    T_n = 
    \sqrt{n}\bigl(\w_p^p(P_n,Q)  - \mathbb  E\bigl(\w_p^p(P_n,Q)\bigr)\bigr) +  \sqrt{n}\bigl(\mathbb E\bigl(\w_p^p(P_n,Q)\bigr)-\w_p^p(P,Q) \bigr)  \ ,
\end{align}
which together with the weak convergence allow us to conclude. 
The case \(p=2\) follows similarly, using that
\begin{equation}\label{eq:decomposition_quantile}
 T_n =  \sqrt{n} \Bigl(  \int_0^1 \bigl( F_n^{-1}(t)^2 - F^{-1}(t)^2 \bigr) dt \Bigr)- 2   \Bigl( \int_0^1 G^{-1}(t)\sqrt{n}\bigl( F_n^{-1}(t) - F^{-1}(t) \bigr) dt\Bigr) \ .
\end{equation}
Since $\int_0^1F_n^{-1}(t)^2dt = (1/n)\sum_{i=1}^n X_i^2$, the first term is uniformly integrable and has zero expectation. The second term  weakly converges to a centered Gaussian random variable by \Cref{lemma:csorgo}. Observe that, in this case, additional moment assumptions are needed only for \(Q\). The same uniform integrability argument allows us to conclude.\\

\noindent
\textit{Two-sample setting:} The proof follows from the same reasoning, with  $T_{n,m}  =  \sqrt{\frac{nm}{n+m}} \bigl(\w_p^p(P_n,Q_m) - \w_p^p(P,Q)\bigr)$,
\begin{align*}
     &L_{n,m} =  \int_0^1 h'_p\bigl( F^{-1}(t) - G^{-1}(t)\bigr)\Bigr[  \sqrt{\frac{m}{n+m}}  \sqrt{n}\bigl( F_n^{-1}(t) - F^{-1}(t) \bigr) - \sqrt{\frac{n}{n+m}}   \sqrt{m}\bigl( G_m^{-1}(t) - G^{-1}(t) \bigr)\Bigr] dt   \\
     &U_{n,m}  =  \int_0^1 h'_p\bigl( F_n^{-1}(t) - G_m^{-1}(t)\bigr)\Bigr[  \sqrt{\frac{m}{n+m}}  \sqrt{n}\bigl( F_n^{-1}(t) - F^{-1}(t) \bigr) - \sqrt{\frac{n}{n+m}}   \sqrt{m}\bigl( G_m^{-1}(t) - G^{-1}(t) \bigr)\Bigr] dt   \ .
\end{align*}
As before, $L_{n,m}\leq T_{n,m}\leq U_{n,m}$ and $U_{n,m}-L_{n,m}\rightarrow_p 0$, whereas $L_{n,m}$ weakly converges by \Cref{lemma:csorgo} to a sum of two 
independent centered Gaussian variables, which is again a centered Gaussian variable. This allows us to conclude using uniform integrability as before. Again, the case $p=2$ follows from the same ideas, but expanding the squares.\\

\noindent
\textbf{Extension to the sliced setting : \ }
We aim to prove \eqref{eq:bias_conv_0}, which can be rewritten as
$$
    \int_{\mathbb S^{d-1}}  \sqrt{n} \Bigl( \mathbb E \left( \w_p^p(P_n^\theta, Q^\theta)\right) -  \w_p^p(P^\theta, Q^\theta) \Bigr)d\sigma(\theta) \rightarrow 0 \ .
$$
Given that $\operatorname{SJ}_{\alpha}(P)<\infty$, then $\operatorname{J}_{\alpha}(P^{\theta})<\infty$ for almost every $\theta \in \mathbb S^{d-1}$. In addition, the moment assumptions on $P,Q$ of Theorem \ref{teor:bias_improved_csorgo_sliced} imply the same moment assumptions for $P^\theta,Q^\theta$  for each $\theta \in \mathbb S^{d-1}$. Thus, from the one-dimensional case, it follows that  
$
        \sqrt{n} ( \mathbb E(\w_p^p(P_n^\theta,Q^\theta)) - \w_p^p(P^\theta,Q^\theta)) \rightarrow 0,
$ 
for almost every $\theta \in \mathbb S^{d-1}$. We have pointwise convergence towards 0 of the function within the integral. Therefore, we can conclude by the dominated convergence theorem if we bound the function by an integrable expression. Define $h_{\theta,n}(t) = 2G_\theta^{-1}(t)$ if $p=2$ and $h_{\theta,n}(t) = p|F_{\theta,n}^{-1}(t)-G_\theta^{-1}(t)|^{p-1}$ if $p\neq 2$.
For each $p$, reasoning as in the one-dimensional case and using the non-negativity property \eqref{eq:bias_conv_0} and Proposition 5.3 in \cite{bobkov2019onedimensional}, it follows that
\begin{align*}
   0&\leq   \sqrt{n} \Bigl( \mathbb E\bigl(\w_p^p(P_n^\theta,Q^\theta)\bigr) - \w_p^p(P^\theta,Q^\theta)\Bigr) \leq  \mathbb E \Bigl( \int_0^1 \Big| h_{\theta,n}(t) \sqrt{n}(F_{\theta,n}^{-1}(t)-F_\theta^{-1}(t))\Big| dt  \Bigr)  \\
    & \leq 5\alpha \ \mathbb E \Bigl(  \int |h_{\theta,n}(t)|^{\beta}dt \Bigr)^{1/\beta} \operatorname{J}_{\alpha}^{1/\alpha} (P^\theta) \ .
\end{align*}
In any case, under the given moment assumptions, there exists a constant  $M<\infty$ such that $ \mathbb E \bigl(  \int |h_{\theta,n}(t)|^{\beta}dt \bigr)^{1/\beta} \leq M$.
Since we are assuming that $\theta \mapsto \operatorname{J}_\alpha(P^\theta) \in L^1(\sigma)$, and $\alpha \geq 1$, applying Jensen's inequality, it follows that  $\theta \mapsto \operatorname{J}_2(P^\theta)^{1/\alpha} \in L^1(\sigma)$, and we can conclude applying dominated convergence theorem.\\

\noindent
\textit{Two-sample setting:} The one-dimensional case ensures that
$ \sqrt{\frac{nm}{n+m}}  \Bigl( \mathbb E\bigl(\w_p^p(P_n^\theta,Q_m^\theta)\bigr) - \w_p^p(P^\theta,Q^\theta)\Bigr) \rightarrow 0  \ \forall \ \theta \in \mathbb S^{d-1}.$
To bound this term by an integrable function, we can decompose the problem as follows
\begin{align*}
   0&\leq    \sqrt{\frac{nm}{n+m}} \Bigl( \mathbb E\bigl(\w_p^p(P_n^\theta,Q_m^\theta)\bigr) - \w_p^p(P^\theta,Q^\theta)\Bigr) \\ 
   &\leq     \sqrt{m} \Bigl( \mathbb E\bigl(\w_p^p(P_n^\theta,Q_m^\theta)\bigr) - \w_p^p(P_n^\theta,Q^\theta)\Bigr) + \sqrt{n} \Bigl( \mathbb E\bigl(\w_p^p(P_n^\theta,Q^\theta)\bigr) - \w_p^p(P^\theta,Q^\theta)\Bigr)
\end{align*}
By the one-sample setting, there exists $C<\infty$ such that the last term is bounded by $CJ _\alpha^{1/\alpha}(P^\theta)$. Similarly, 
\begin{align*}
     \sqrt{n} \Bigl( \mathbb E\bigl(\w_p^p(P_n^\theta,Q_m^\theta)\bigr) - \w_p^p(P_n^\theta,Q^\theta)\Bigr)     
    &\leq  5\alpha' \ \mathbb E \Bigl(  \int |h_{\theta,n,m}(t)|^{\beta'}dt \Bigr)^{1/\beta'} \operatorname{J}_{\alpha'}^{1/\alpha'} (Q^\theta) = C' \operatorname{J}_{\alpha'}^{1/\alpha'} (Q^\theta)\ ,
\end{align*}
where $h_{\theta,n,m}(t) = 2F_{\theta,n}^{-1}(t)$ if $p=2$ and $h_{\theta,n}(t) = p|F_{\theta,n}^{-1}(t)-G_{\theta,m}^{-1}(t)|^{p-1}$ if $p\neq 2$, and, therefore, $C'<\infty$. Since $\theta \mapsto CJ _\alpha^{1/\alpha}(P^\theta)+C'J _{\alpha'}^{1/\alpha'}(Q^\theta)$ is integrable under the given assumptions, we can conclude as before. \qed

\section{Proofs of Section \ref{sect:slicing}}\label{appendixC}

\subsection{Proof of Proposition \ref{teor:finite_slicing}}

The proof is inspired by the bootstrap results of \cite{bickel1981someAsymptotic}.
First, given a random variable $\Theta$ drawn from $\sigma$, define $Z_n(\Theta) = \w_p^p(P_n^\Theta,Q^\Theta)$ and $Z(\Theta) = \w_p^p(P^\Theta,Q^\Theta)$. Throughout this proof, $\mathbb{E}_{\Theta}$ denotes expectation taken only over $\Theta_1, \ldots, \Theta_k$, whereas $\mathbb{E}$ represents  expectation over all random variables. With this convention,
\begin{align}
        S_{n,k} = \frac{1}{\sqrt k} \sum_{i=1}^k \Bigl( Z_n(\Theta_i) - \mathbb E_{\Theta}(Z_n(\Theta_i))\Bigr) \ , \quad \textnormal{ and let } \quad    S_{k} = \frac{1}{\sqrt k} \sum_{i=1}^k \Bigl( Z(\Theta_i) - \mathbb E_{\Theta}(Z(\Theta_i))\Bigr) \ . 
        \label{eq:S_{nk}} 
\end{align}
By definition, $ \mathbb E(Z(\Theta))= \sw_p^p(P,Q)$ and $\operatorname{Var}(Z(\Theta)) = \int_{\mathbb S^{d-1}} W_p^{2p}(P^\theta,Q^\theta) d\sigma(\theta) - \sw_p^{2p}(P,Q) $, 
which is finite, since $\w_p^{2p}(P^\theta,Q^\theta)\leq \w_p^{2p}(P,Q) \leq \w_{2p}^{2p}(P,Q)$ for every $\theta\in\mathbb S^{d-1}$. Therefore, a straightforward application of the CLT yields  $S_k\rightsquigarrow N(0,w_{P,Q}^2)$ and second-order moment convergence, which implies  $W_2(\mathcal L(S_k),N(0,w_{P,Q}^2) )\rightarrow 0$, by Lemma 8.3 in \cite{bickel1981someAsymptotic}. Consider now the coupling given by independent variables $X_1\ldots,X_n$ and $\Theta_1,\ldots,\Theta_k$.     
\begin{align}
        \w_2^2( \mathcal L(S_{n,k}),\mathcal L(S_k)) &\leq \mathbb E \left(  \Bigl(\frac{1}{ \sqrt k} \sum_{i=1}^k Z_n(\Theta_i) - \mathbb E_{\Theta}(Z_n(\Theta_i)) -  \bigl(Z(\Theta_i) - \mathbb E_{\Theta}(Z(\Theta_i))\bigr)\Bigr) ^2\right)\notag\\
        &= \frac{1}{k} \sum_{i=1}^k \mathbb E\left( \bigl( Z_n(\Theta_i)  -  Z(\Theta_i) -\mathbb E_{\Theta}(Z_n(\Theta_i)-Z(\Theta_i))\bigr)^2\right)  \notag\\
       & \leq  \mathbb E\left( \bigl( Z_n(\Theta)  -  Z(\Theta)\bigr)^2\right)  = \mathbb  E\left(  \bigl( \w_p^p(P_n^\Theta,Q^\Theta)-\w_p^p(P^\Theta,Q^\Theta)\bigr)^2 \right) \label{eq:slicing_last_ineq}
         \end{align}
where the expected values of the cross terms in the equality vanish due to the independence of $\Theta_i$ and $\Theta_j$ for $i\neq j$.
From the bound \eqref{eq:convex_basic_ineq} and the triangle inequality, for every $\theta\in\mathbb S^{d-1}$,
\begin{align*}
    \big| \w_p^p(P_n^\theta,Q^\theta)-\w_p^p(P^\theta,Q^\theta) \big| 
    & \leq p \bigl( \w_p^{p-1}(P_n^\theta,Q^\theta)+\w_p^{p-1}(P^\theta,Q^\theta)\Bigr)   \w_p(P_n^\theta,P^\theta) \ .
\end{align*}
By the above inequality, \eqref{eq:elemental_ineq_p} and  Hölder's inequality with $q_1=\frac{p}{p-1}$, $q_2=p$, \eqref{eq:slicing_last_ineq} is bounded, up to a constant, by 
\begin{align*}
    \Bigl(\mathbb  E\bigl( & \w_p^{2p}(P_n^\Theta,Q^\Theta)\bigr) + \mathbb E\bigl(  \w_p^{2p}(P^\Theta,Q^\Theta)\bigr)\Bigr)^{\frac{p-1}{p}} \mathbb  E\bigl(  \w_p^{2p}(P_n^\Theta,P^\Theta)\bigr)^{\frac{1}{p}} 
\end{align*}
Using again that $\w_p^{2p}\leq \w_{2p}^{2p}$, and since we are assuming that $P,Q$ have finite moment of order $2p$, the first term is bounded, whereas the second term converges to zero by  Lemma \ref{lemma:convergence_EW_p}. \qed

\subsection{Proof of Corollary \ref{coro_slicing}}

 Let $S_{n,k}$ and $S_k$ be as in \eqref{eq:S_{nk}}, and define
 $
    R_n = \sqrt n \Bigl( \sw_p^p(P_n,Q)-\mathbb E \bigl(\sw_p^p(P_n,Q)\bigr) \Bigr) 
$. Then, the left-hand side of \eqref{TCL:finite_sliced_one_sample} can be decomposed as
\begin{align*}
  \sqrt{\frac{n}{k+n}} S_{n,k} +  \sqrt{\frac{k}{k+n}}R_n +  \sqrt{\frac{k}{k+n}} \sqrt n \Bigl(\mathbb E \bigl(\sw_p^p(P_n,Q)\bigr)-\sw_p^p(P,Q)\Bigr) \ .
\end{align*}
Under the given assumptions, Theorem \ref{teor:bias_improved_csorgo_sliced} ensures that the last term converges to zero, while Proposition \ref{teor:finite_slicing} establishes weak convergence of the first term, and Theorem \ref{teor:CLT} guarantees weak convergence of the second term. Therefore, the only thing left to prove is that the sum of both sequences weakly converges to the sum of two independent Gaussian random variables with the distribution of the previous weak limits. This convergence is not straightforward since the first and second terms are not independent. 
The variable $S_{n,k}$ depends both on $X_1,\ldots,X_n$ and $\Theta_1,\ldots,\Theta_k$, whereas $S_k$ only depends on  $\Theta_1,\ldots,\Theta_k$ and $R_n$ on $X_1,\ldots,X_k$. Therefore, 
\begin{align*}
  &  W_2\Biggl( \mathcal L\Bigl( \sqrt{\frac{n}{k+n}} S_{n,k} + \sqrt{\frac{k}{k+n}} R_{n} \Bigr),  \sqrt{1-\tau}\cdot N(0,w_{P,Q}^2) +  \sqrt{\tau}\cdot N(0,v_{P,Q}^2) \Biggr) \\
   & \leq   W_2\Biggl( \mathcal L\Bigl( \sqrt{\frac{n}{k+n}} S_{n,k} + \sqrt{\frac{k}{k+n}} R_n \Bigr),\mathcal L\Bigl( \sqrt{\frac{n}{k+n}} S_{k} + \sqrt{\frac{k}{k+n}} R_n \Bigr) \Biggr)\\
   & +   W_2\Biggl( \mathcal L\Bigl( \sqrt{\frac{n}{k+n}} S_{k} + \sqrt{\frac{k}{k+n}} R_n \Bigr), \sqrt{1-\tau}\cdot N(0,w_{P,Q}^2) +  \sqrt{\tau}\cdot N(0,v_{P,Q}^2) \Biggr)
\end{align*}
The first term is bounded by $ \w_2( \mathcal L(S_{n,k}),\mathcal L(S_k))$, which converges to zero as shown in the proof of Proposition \ref{teor:finite_slicing}. For the second term, Proposition \ref{teor:finite_slicing} ensures that $W_2(\mathcal L(S_k),N(0,w_{P,Q}^2) )\rightarrow 0$, while Theorem \ref{teor:CLT} demonstrates that $W_2(\mathcal L(R_n),N(0,v_{P,Q}^2) )\rightarrow 0$. Since $S_k$ and $R_n$ are independent, it is easy to conclude 
considering the coupling $(S_k,U,R_n,V)$, where $\mathcal L(U)= N(0,w_{P,Q}^2)$,  $\mathcal L(V)= N(0,v_{p,Q}^2)$, $(S_k,U)$ is the optimal coupling for $\mathcal L(S_k),N(0,w_{P,Q}^2)$,  $(R_n,V)$  optimal coupling for $\mathcal L(S_k),N(0,w_{P,Q}^2)$, and $(S_k,U)$ and $(R_n,V)$ are independent. \qed

\section{Proofs of Section \ref{sect:variance}}\label{appendixD}

\subsection{Proof of Proposition \ref{coro:consistW}}

The proof proceeds by showing that each term in the following decomposition converges in probability to zero.
\begin{align}
    \hat w^2_{P_n,Q} - w^2_{P,Q} &= \frac{1}{k}\sum_{i=1}^k \w_p^{2p}(P_n^{\Theta_i},Q^{\Theta_i}) - \frac{1}{k}\sum_{i=1}^k \w_p^{2p}(P^{\Theta_i},Q^{\Theta_i}) \label{eq:proof_6.1_1}\\
    &+ \frac{1}{k}\sum_{i=1}^k \w_p^{2p}(P^{\Theta_i},Q^{\Theta_i}) - \int_{\mathbb S^{d-1}} \w_p^{2p}(P^{\theta},Q^{\theta})d\sigma(\theta) \label{eq:proof_6.1_2}\\
    &-\Bigl(\frac{1}{k}\sum_{i=1}^k \w_p^{p}(P_n^{\Theta_i},Q^{\Theta_i}) \Bigr)^2+ \Bigl( \int_{\mathbb S^{d-1}} \w_p^p(P^{\theta},Q^{\theta})d\sigma(\theta)  \Bigr)^2\label{eq:proof_6.1_3}
\end{align}
For \eqref{eq:proof_6.1_2}, convergence follows from the strong law of large numbers. For \eqref{eq:proof_6.1_1}, it follows from Markov's inequality and 
\begin{align*}
    \mathbb E \Big | \frac{1}{k}&\sum_{i=1}^k \w_p^{2p}(P_n^{\Theta_i},Q^{\Theta_i}) - \frac{1}{k}\sum_{i=1}^k \w_p^{2p}(P^{\Theta_i},Q^{\Theta_i}) \Big| \leq  \mathbb E  \Bigl( \int_{\mathbb S^{d-1}} \Big| \w_p^{2p}(P_n^{\theta},Q^{\theta}) - \w_p^{2p}(P^{\theta},Q^{\theta}) \Big|d\sigma(\theta) \Bigr) \rightarrow 0 \ ,
\end{align*}
where the convergence follows as in the proof of Proposition \ref{teor:finite_slicing}. Similarly, we can prove
$
\frac{1}{k}\sum_{i=1}^k \w_p^p(P_n^{\Theta_i},Q^{\Theta_i}) - \frac{1}{k}\sum_{i=1}^k \w_p^{p}(P^{\Theta_i},Q^{\Theta_i}) \rightarrow_p 0 \ ,
$
and, by the strong law of large numbers, we know that 
$
 \frac{1}{k}\sum_{i=1}^k \w_p^{p}(P^{\Theta_i},Q^{\Theta_i}) \rightarrow_{a.s.} \int_{\mathbb S^{d-1}} \w_p^p(P^{\theta},Q^{\theta})d\sigma(\theta)  
$. Slutski's lemma allow us to conclude convergence in probability to zero of \eqref{eq:proof_6.1_3}. 
\qed

\subsection{Proof of Proposition \ref{prop:consistencia}}

For each $\theta\in\mathbb S^{d-1}$, given any $c$-concave optimal transport potentials $\phi_n^\theta,\phi^\theta$,  let $\phi_{n,X_1}^\theta$, $\phi_{X_1}^\theta$ be the random optimal transport potentials fixed at $X_1$ defined as in \Cref{prop:regularidad_potenciales}, which are well-defined $\mathbb P$-almost surely. As in previous proofs, we can assume without loss of generality that are well-defined everywhere. In addition, denote by $\phi_{x_0}^\theta$ the optimal transport potential fixed at $x_0\in\operatorname{int}(\supp{P})$ defined by \eqref{defn:fixed_potentials_x0}. By Corollary \ref{corollary:integrability}, we can assume \eqref{eq:int_potentials_limit_x0} for $q=2(1+\gamma)$.
Since the covariance is invariant by the addition of a constant, for each $\theta,\eta\in\mathbb S^{d-1}$, $ \operatorname{Cov}_{P_n}(\phi^{\theta}_{n},\phi^{\eta}_{n}) = \operatorname{Cov}_{P_n}(\phi^{\theta}_{n,X_1},\phi^{\eta}_{n,X_1})$, $\operatorname{Cov}_{P_n}(\phi^{\theta}_{X_1},\phi^{\eta}_{X_1})  = \operatorname{Cov}_{P_n}(\phi^{\theta}_{x_0},\phi^{\eta}_{x_0})$ and $\operatorname{Cov}_{P}(\phi^{\theta}_{x_0},\phi^{\eta}_{x_0})  = \operatorname{Cov}_{P}(\phi^{\theta},\phi^{\eta})$, where in the last two equalities we also use the uniqueness of the limit potentials up to a constant of Corollary 2.7 in \cite{delBarrio2024CLTgeneral}.
Throughout the remainder of the proof, given any optimal transport potential 
$\phi^\theta$, we write with a slight abuse, $\phi^\theta(x)= \phi^\theta(\langle\theta,x\rangle)$. With this convention,
    \begin{align*}
        \hat v_{P_n,Q}^2-v_{P,Q}^2 
        &=\frac{1}{k^2} \sum_{i,j=1}^k \operatorname{Cov}_{P_n}(\phi^{\Theta_i}_{n,X_1},\phi^{\Theta_j}_{n,X_1})- \frac{1}{k^2} \sum_{i,j=1}^k \operatorname{Cov}_{P_n}(\phi^{\Theta_i}_{X_1},\phi^{\Theta_j}_{X_1}) + \\
        &  + \frac{1}{k^2} \sum_{i,j=1}^k \operatorname{Cov}_{P_n}(\phi^{\Theta_i}_{x_0},\phi^{\Theta_j}_{x_0})- \int_{\mathbb S^{d-1}} \int_{\mathbb S^{d-1}} \operatorname{Cov}_{P}(\phi^\theta_{x_0},\phi^\eta_{x_0}) d\sigma(\theta) d\sigma(\eta)  \\
         &=\frac{1}{k^2n} \sum_{i,j=1}^k\sum_{l=1}^n \Bigl(\phi^{\Theta_i}_{n,X_1}(X_l) \phi^{\Theta_j}_{n,X_1}(X_l) - \phi^{\Theta_i}_{X_1}(X_l) \phi^{\Theta_j}_{X_1}(X_l) \Bigr) 
         \tag{A} \\
         &  - \Biggl( \Bigl(\frac{1}{kn} \sum_{i=1}^k\sum_{l=1}^n  \phi^{\Theta_i}_{n,X_1}(X_l)\Bigr)^2 - \Bigl(\frac{1}{kn} \sum_{i=1}^k \sum_{l=1}^n \phi^{\Theta_i}_{X_1}(X_l)\Bigr)^2 \Biggr)  \tag{B}\\
            &  + \frac{1}{k^2n} \sum_{i,j=1}^k\sum_{l=1}^n\phi^{\Theta_i}_{x_0}(X_l) \phi^{\Theta_j}_{x_0}(X_l)  - \int_{\mathbb S^{d-1}} \int_{\mathbb S^{d-1}} \int_{\mathbb R^{d}} \Bigl(\phi_{x_0}^\theta(x)\phi_{x_0}^\eta(x)\Bigr) \ dP(x) d\sigma(\theta) d\sigma(\eta) \tag{C}\\
              &   - \Biggl(  \Bigl(\frac{1}{kn} \sum_{i=1}^k\sum_{l=1}^n  \phi^{\Theta_i}_{x_0}(X_l)\Bigr)^2   -  \Bigl(\int_{\mathbb S^{d-1}} \int_{\mathbb R^{d}} \phi^\theta_{x_0}(x)\ dP(x) d\sigma(\theta) \Bigr)^2 \Biggr)\tag{D}
    \end{align*}
    We will prove convergence in probability to zero of each term in the decomposition. First, we  start with (D). From the continuous mapping theorem, it suffices to see that 
  $
        \frac{1}{kn} \sum_{i=1}^k \sum_{l=1}^n \phi^{\Theta_i}_{x_0}(X_l)   \rightarrow_p \int_{\mathbb S^{d-1}} \int_{\mathbb R^{d}} \phi^\theta_{x_0}(x)\ dP(x) d\sigma(\theta) \ .
$
Consider the function $h(\theta,x)=\phi^\theta_{x_0}(x)$. Since $\Theta_1,\ldots,\Theta_k$ and $X_1,\ldots,X_n$ are independent sequences, the problem can be reformulated in terms of two sample U-statistics
(also known as generalized U-statistics). 
The above convergence can be rewritten as 
      \begin{align*}
        \frac{1}{kn} \sum_{i=1}^k \sum_{l=1}^n h(\Theta_i,X_l)   \rightarrow_p \mathbb E\bigl( h(\Theta_1,X_1)\bigr) \ .
    \end{align*}
    In fact, existing versions of the strong law of large numbers for generalized U-statistics allow us to prove convergence almost surely. If $n$ and $k$ are of the same order, i.e. $n/(n+k)\rightarrow \tau\in (0,1)$, Theorem 3.1.2. in \cite{Korolyuk2013theoryUstatistic} proves a.s. convergence provided that  $\mathbb E\bigl( |h(\Theta_1,X_1)|\bigr)<\infty$. Assuming only $n,k\rightarrow \infty$, \cite{Sen1977Ustatistics} proves a.s. convergence under the stronger assumption  $\mathbb E\bigl( |h(\Theta_1,X_1)|\log_+(|h(\Theta_1,X_1)|)\bigr)<\infty$. Both conditions are verified, since
    \begin{equation}\label{eq:boundUstat}
        \mathbb E\bigl( |h(\Theta_1,X_1)|^{2(1+\gamma)}\bigr) =   \| \phi_{x_0}^\theta \|_{L^{2(1+\gamma)}(\sigma\times P)}^{2(1+\gamma)}<\infty \ .
    \end{equation}
    For the term (C), consider the following decomposition
    \begin{align*}
        &\frac{1}{k^2n} \sum_{i,j=1}^k\sum_{l=1}^n\phi^{\Theta_i}_{x_0}(X_l) \phi^{\Theta_j}_{x_0}(X_l) =  \frac{1}{k}\Biggl( \frac{1}{kn} \sum_{i=1}^k\sum_{l=1}^n\phi^{\Theta_i}_{x_0}(X_l)^2  \Biggr)
         +   \frac{2\binom{k}{2}}{k^2} \Biggl(  \frac{1}{\binom{n}{1}\binom{k}{2}}\sum_{1\leq i<j\leq k}\sum_{l=1}^n\phi^{\Theta_i}_{x_0}(X_l) \phi^{\Theta_j}_{x_0}(X_l) \Biggr)
    \end{align*}
    The bound \eqref{eq:boundUstat}, and the same reasoning as before, allow us to conclude that the term inside the first parenthesis converges a.s. to $\mathbb E\bigl( h(\Theta_1,X_1)^{2}\bigr)$. Since it is multiplied by a factor $1/k\rightarrow 0$, the first term converges a.s. to zero. Now, given that $2\binom{k}{2}/k^2 = (k-1)/k \rightarrow 1$, if we denote $g(\theta,\theta',x) = \phi^{\theta}_{x_0}(x) \phi^{\theta'}_{x_0}(x) $, 
    it suffices to prove 
     \begin{equation*}
         \frac{1}{\binom{n}{1}\binom{k}{2}}\sum_{1\leq i<j\leq k}\sum_{l=1}^n g(\Theta_i,\Theta_j,X_l)  \rightarrow_p \mathbb E\bigl( g(\Theta_1,\Theta_2,X_1)\bigr)
     \end{equation*}
    to conclude (C). Given that $g(\theta,\theta',x)$ is symmetric with respect to $\theta,\theta'$, we can apply again strong law of large numbers for two sample U-statistics to conclude almost sure convergence, since
    \begin{align*}
        \mathbb E \Bigl( |g(\Theta_1,\Theta_2,X_1)|^{1+\gamma} \Bigr) &\leq          \mathbb E \Bigl( |\phi_{x_0}^{\Theta_1}(X_1)|^{2(1+\gamma)} \Bigr)^{1/2}\mathbb E \Bigl( |\phi_{x_0}^{\Theta_2}(X_1)|^{2(1+\gamma)} \Bigr)^{1/2} = \| \phi_{x_0}^\theta \|_{L^{2(1+\gamma)}(\sigma\times P)}^{2(1+\gamma)}<\infty  \ . 
    \end{align*}
Next, to show the convergence of (B), we will prove the convergence in mean  
\begin{align} \label{eq:conv_mean1}
    \mathbb E \Bigl|\frac{1}{kn} \sum_{i=1}^k\sum_{l=1}^n  \phi^{\Theta_i}_{n,X_1}(X_l) - \frac{1}{kn} \sum_{i=1}^k\sum_{l=1}^n  \phi^{\Theta_i}_{X_1}(X_l)\Bigr| \rightarrow 0  \ .
\end{align}
By Markov's inequality, \eqref{eq:conv_mean1} implies convergence in probability to zero. Using Proposition \ref{prop:integrability_potentials}, it is easy to see that the sequence $\frac{1}{kn} \sum_{i=1}^k\sum_{l=1}^n  \phi^{\Theta_i}_{n,X_1}(X_l) + \frac{1}{kn} \sum_{i=1}^k\sum_{l=1}^n  \phi^{\Theta_i}_{X_1}(X_l)$ is bounded in probability, and therefore, (B) follows multiplying.
To prove \eqref{eq:conv_mean1}, note that 
\begin{align*}
   \phi_{n,X_1}^\theta(X_1) = \phi_{X_1}^\theta(X_1)  = 0 \ , \quad \phi_{n,X_1}^\theta(X_2)\underset{d}{=}\ldots\underset{d}{=}\phi_{n,X_1}^\theta(X_n) \quad \textnormal{ and }\quad  \phi_{X_1}^\theta(X_2)\underset{d}{=}\ldots\underset{d}{=}\phi_{X_1}^\theta(X_n) \ .
\end{align*}
Hence, \eqref{eq:conv_mean1} can be bounded by 
\begin{align*} 
    \frac{1}{kn} \mathbb E \Bigl(\sum_{i=1}^k\sum_{l=1}^n  | \phi^{\Theta_i}_{n,X_1}(X_l) - \phi^{\Theta_i}_{X_1}(X_l)| \Bigr)& = \frac{1}{n} \mathbb E \Bigl(\sum_{l=1}^n  | \phi^{\Theta_1}_{n,X_1}(X_l) - \phi^{\Theta_1}_{X_1}(X_l)| \Bigr) =  \frac{n-1}{n} \mathbb E \Bigl( \big|\phi_{n,X_1}^{\Theta_1}(X_2)-\phi_{X_1}^{\Theta_1}(X_2) \big|\Bigr)
\end{align*}
By property (4) in Proposition \ref{prop:regularidad_potenciales}, we know that $|\phi_{n,X_1}^{\Theta_1}(X_2) - \phi_{X_1}^{\Theta_1}(X_2)|\rightarrow 0$ $\mathbb P$-a.s. Uniform integrability follows from Proposition \ref{prop:integrability_potentials}, since
\begin{align}\label{eq:reasoning_uniform_int}
\mathbb E  \Bigl( \big|\phi_{n,X_1}^{\Theta_1}(X_2)-\phi_{X_1}^{\Theta_1}(X_2) \big|^{1+\gamma}\Bigr)    &= \mathbb E\Bigl(\int_{\mathbb S^{d-1}} \frac{1}{n-1}\sum_{l=2}^n \big| \phi^\theta_{n,X_1}(X_l) - \phi^\theta_{X_1}(X_l) \Big|^{1+\gamma} d\sigma(\theta)\Bigr)\notag\\
   & \leq C_{1+\gamma} \Bigl( \frac{n}{n-1}\|\phi^\theta_{n,X_1}\|_{L^{1+\gamma}(\sigma\times P_n)}^{1+\gamma} + \| \phi^\theta_{X_1}\|_{L^{1+\gamma}(\sigma\times P)}^{1+\gamma}\Bigr)
\end{align}
is uniformly bounded. A similar conclusion is reached for $(A)$. To prove convergence in mean, 
\begin{align*}
    &\mathbb E \Big| \frac{1}{k^2n} \sum_{i,j=1}^{k}\sum_{l=1}^n \Bigl(\phi^{\Theta_i}_{n,X_1}(X_l) \phi^{\Theta_j}_{n,X_1}(X_l) - \phi^{\Theta_i}_{X_1}(X_l) \phi^{\Theta_j}_{X_1}(X_l) \Bigr)  \Big| \\
    &  \leq \frac{2}{k^2n} \sum_{1\leq i<j\leq n}^k\sum_{l=1}^n \mathbb E \big| \phi^{\Theta_i}_{n,X_1}(X_l) \phi^{\Theta_j}_{n,X_1}(X_l) - \phi^{\Theta_i}_{X_1}(X_l) \phi^{\Theta_j}_{X_1}(X_l)  \big|  +  \frac{1}{k} \frac{1}{kn} \sum_{i=1}^k\sum_{l=1}^n \mathbb E \big|\phi^{\Theta_i}_{n,X_1}(X_l)^2 - \phi^{\Theta_i}_{X_1}(X_l)^2  \big| \\
     &= \frac{(k-1)(n-1)}{kn}\mathbb E \big|\phi^{\Theta_1}_{n,X_1}(X_2) \phi^{\Theta_2}_{n,X_1}(X_2) - \phi^{\Theta_1}_{X_1}(X_2) \phi^{\Theta_2}_{X_1}(X_2)  \big| +  \frac{1}{k} \frac{n-1}{n}  \mathbb E \big|\phi^{\Theta_1}_{n,X_1}(X_2)^2 - \phi^{\Theta_1}_{X_1}(X_2)^2  \big| 
    \end{align*}
By Proposition \ref{prop:regularidad_potenciales}, $\phi^{\Theta_1}_{n,X_1}(X_2) - \phi^{\Theta_1}_{X_1}(X_2)$ $\mathbb P$-a.s. Therefore, 
we have pointwise convergence of both sequences inside the expectations. Uniform integrability can be shown as in \eqref{eq:reasoning_uniform_int}, using again \Cref{prop:integrability_potentials} and a further application of H\"older's inequality.\\

\noindent
\textit{Two-sample setting:} By symmetry, it suffices to see $\hat v^2_{P_n,Q_m}\rightarrow_p v^2_{P,Q}$. The proof follows from the same arguments, since the uniform convergence and uniform integrability properties in the above proof also apply to the optimal transport potentials $\phi_{n,m,X_1}^\theta$, defined as in \Cref{prop:regularidad_potenciales}.
\qed

\section{Efficient estimation of the asymptotic variance}\label{appendix:efficient}

In this section, we detail the computations required to obtain the value of $\sw_{p,k}^p(P_n,Q_m)$ and the variance estimators  $\hat w_{P_n,Q_m}^2$, $\hat v_{P_n,Q_m}^2$ and $\hat v_{Q_m,P_n}^2$ in the two-sample setting. We restrict our attention to the prominent case $p=2$, where the $c$-concave optimal transport potentials can be efficiently estimated using the relationship between $c$-concavity and Legendre duality. For general $p > 1$, different strategies may be adopted to find the $c$-concave optimal transport potentials that solve the dual problem. However, it remains unclear whether a straightforward solution exists that can be efficiently computed based on the one-dimensional optimal transport plan, as in the case of $p=2$.\\

\noindent
Let $x_1,\ldots,x_n$, $y_1,\ldots,y_m \in \mathbb R^d$ and $\theta_1,\ldots,\theta_k\in\mathbb S^{d-1}$ be the given samples. Denote the projected values by 
$
    s_i^l = \langle \theta_l,x_i\rangle , \  
    t_j^l =\langle \theta_l,y_j\rangle, \  \textnormal{for } i\in [n],  \ j\in[m], \ l\in[k]  ,
$
and let $\{s_{(i)}^l\}_{i=1}^n$ and $\{t_{(j)}^l\}_{j=1}^m$  be the values of the ordered statistic. For each $l\in[k]$ we need to compute 
\begin{align*}    \w_2^2(P_{n}^{\theta_l},Q_{m}^{\theta_l}) & = \int_0^1 \Bigl(F_{\theta_l,n}^{-1}(u)-G_{\theta_l,m}^{-1}(u)\Bigr)^2du= \sum_{i=1}^{n} \sum_{j=1}^{m} R_{i,j}(s_{(i)}^l-t_{(j)}^l )^2  \ ,
\end{align*}
where $ R_{i,j}= \ell_1 \bigl((\frac{i-1}{n},\frac{i}{n}]\cap(\frac{j-1}{m},\frac{j}{m}]\bigr)$. The matrix $R=\{R_{i,j}\}_{i,j}$ only depends on $n$ and $m$, it must be computed just once at the beginning. Once we have the values of $\{\w_2^2(P_{n}^{\theta_l},Q_{m}^{\theta_l})\}_{l=1}^k$, we can compute $\sw_{2,k}^2(P_n,Q_m)$ as in  \eqref{eq:defn_finite_sliced},  and 
$
 \hat w_{P_n,Q_m}^2 =   \frac{1}{k}\sum_{l=1}^k \w_2^4(P_{n}^{\theta_l},Q_{m}^{\theta_l}) - \sw_{p,k}^{2p}(P_n,Q_m) \ .
$
To compute $\hat v_{P_n,Q_m}^2$, we can reuse part of the previous computations. Our approach is based on the characterization of optimal transport plans in terms of lower semicontinuous convex functions for the cost $c=c_2$ (see \cite{villani2003Topics} for details). In this case, the $c$-concave optimal transport potentials from $P_n^{\theta_l}$ to $Q_m^{\theta_l}$ can be expressed as $\phi_{n,m}^{\theta_l}(x)= x^2- 2\varphi_{n,m}^{\theta_l}(x)$, where $\varphi_{n,m}^{\theta_l}$ is convex and lower semicontinuous. Given an optimal transport plan $\pi_{n,m}^{\theta_l}$, if $\supp{\pi_{n,m}^{\theta_l}} \subset \partial^c \phi_{n,m}^{\theta_l} = \partial \varphi_{n,m}^{\theta_l}$, where $\partial$ denotes the subgradient, then $\phi_{n,m}^{\theta_l}$ is optimal. Therefore, we need to build $\varphi_{n,m}^{\theta_l}$, convex and lower semicontinuous verifying $\supp{\pi_{n,m}^{\theta_l}} \subset  \partial \varphi_{n,m}^{\theta_l}$.
The choice of $\varphi_{n,m}^{\theta_l}$ is not unique. We  consider the only piecewise linear continuous function such that $\varphi_{n,m}^{\theta_l}(s_{(1)}^l)=0$ and 
\begin{align*}
(\varphi_{n,m}^{\theta_l})'(s) =  \left\lbrace\begin{array}{ccc}
    t_{(1)}^l & \textnormal{ if } & s<s_{(1)}^l \\
    t_{(r^l(i))}^l &  \textnormal{ if } & s_{(i)}^l<s<s_{(i+1)}^l \\
    t_{(m)}^l & \textnormal{ if } &  s>s_{(n)}^l 
\end{array}\right. \ ,  \quad  \textnormal{where } \   r(i)= \max \{ j=1,\ldots,m : R_{i,j}>0\} \ .
\end{align*}
To compute $\hat v_{P_n,Q_m}^2$, we are only interested in the values of $\varphi_{n,m}^{\theta_l}$ at the points $\{s_{(i)}^l\}_{i=1}^n$. Since we are imposing $\varphi_{n,m}^{\theta_l}(s_{(1)}^l)=0$, we can iteratively compute 
\begin{align*}          
\varphi_{n,m}^{\theta_l}(s_{(i+1)}^l)=\varphi_{n,m}^{\theta_l}(s_{(i)}^l) +  t_{(r^l(i))}^l ( s_{(i+1)}^l-s_{(i)}^l) \ , \quad i=1,\ldots,n-1.
\end{align*}
Finally, once we have computed all these values, we can obtain the values $\phi_{n,m}^{\theta_l}(s_{(i)}^l)$, which allows us to compute $\hat v_{P_n,Q_m}^2$ as in \Cref{prop:consistencia}. Symmetric computations allow us to obtain $\hat v_{Q_m,P_n}^2$.

\end{document}